\def\version{0.66}
\def\journal{arXiv}
\def\titlep{Classification of sub-Cuntz states}
\documentclass[11pt]{article}
\usepackage{graphicx,ifthen,amssymb,epic,eepic,color}

\newcommand{\qed}{\hbox{\rule[-2pt]{3pt}{6pt}}}
\newcommand{\qedh}{\hfill\qed \\}

\font\germ=eufm10 at12pt

\def\goth#1{\hbox{\germ#1}}

\newcommand{\vv}{\vspace{.3in}}

\newtheorem{Thm}{Theorem}[section]
\newtheorem{fig}[Thm]{Figure}

\newtheorem{rem}[Thm]{Remark}

\newtheorem{ex}[Thm]{Example}
\newtheorem{defi}[Thm]{Definition}
\newtheorem{lem}[Thm]{Lemma}

\newtheorem{prop}[Thm]{Proposition}

\newtheorem{cor}[Thm]{Corollary}
\newtheorem{fact}[Thm]{Fact}

\newcommand{\ww}{\vv\noindent}

\newcommand{\kn}{\Large\bf
$K\hspace{-.4cm} N$
\Large\bf\vv }

\def\cal#1{\mathcal #1}
\def\con{{\cal O}_{n}}
\def\coni{{\cal O}_{\infty}}
\def\edot{=1,\ldots,n}
\def\pr{{\it Proof.}\quad}

\def\co#1{{\cal O}_{#1}}
\def\ltn{\ell^{2}({\Bbb N})}

\def\disp#1{{\displaystyle #1}}
\addtocounter{footnote}{1}
\def\sftt#1{
\setcounter{equation}{0}
\addtocounter{footnote}{1}
\section{#1}
}

\def\ssft#1{\subsection{#1}}
\def\sssft#1{\subsubsection{#1}}

\begin{document}
%
% Personal data
%
\def\autherp{Katsunori Kawamura}
\def\emailp{e-mail: 
%kawamura@kurims.kyoto-u.ac.jp.
kawamurakk3@gmail.com
}
\def\addressp{{\small {\it College of Science and Engineering, 
Ritsumeikan University,}}\\
{\small {\it 1-1-1 Noji Higashi, Kusatsu, Shiga 525-8577, Japan}}
}

\def\scm#1{S({\Bbb C}^{N})^{\otimes #1}}
\def\mqb{\{(M_{i},q_{i},B_{i})\}_{i=1}^{N}}
\newcommand{\mline}{\noindent
\thicklines
\setlength{\unitlength}{.1mm}
\begin{picture}(1000,5)
\put(0,0){\line(1,0){1250}}
\end{picture}
\par
 }
\def\ptimes{\otimes_{\varphi}}
\def\delp{\Delta_{\varphi}}
\def\delf{\Delta_{{\Bbb F}}}
\def\delps{\Delta_{\varphi^{*}}}
\def\gamp{\Gamma_{\varphi}}
\def\gamps{\Gamma_{\varphi^{*}}}
\def\sem{\textsf{M}}
\def\hdelp{\hat{\Delta}_{\varphi}}
\def\tilco#1{\tilde{\co{#1}}}
% Boldfont
\def\ba{\mbox{\boldmath$a$}}
\def\bb{\mbox{\boldmath$b$}}
\def\bc{\mbox{\boldmath$c$}}
\def\be{\mbox{\boldmath$e$}}
\def\bp{\mbox{\boldmath$p$}}
\def\bq{\mbox{\boldmath$q$}}
\def\bu{\mbox{\boldmath$u$}}
\def\bv{\mbox{\boldmath$v$}}
\def\bw{\mbox{\boldmath$w$}}
\def\bx{\mbox{\boldmath$x$}}
\def\by{\mbox{\boldmath$y$}}
\def\bz{\mbox{\boldmath$z$}}

%%%%%%%%%%%%%%%%%%%%%%%%%%%%%%
\def\titlepage{%\vspace{-4cm}

\noindent
{\bf 
\noindent
\thicklines
\setlength{\unitlength}{.1mm}
\begin{picture}(1000,0)(0,-300)
\put(0,0){\kn \knn\, for \journal\, Ver.\,\version}
\put(0,-50){\today,\quad {\rm file:} {\rm {\small \textsf{tit01.txt,\, J1.tex}}}}
\end{picture}
}
\vspace{-2.5cm}
\quad\\
{\small 
\footnote{
\begin{minipage}[t]{6in}
directory: \textsf{\fileplace}, \\
file: \textsf{\incfile},\, from \startdate
\end{minipage}
}}
\quad\\
\framebox{
%\noindent
\begin{tabular}{ll}
\textsf{Title:} &
\begin{minipage}[t]{4in}
\titlep
\end{minipage}
\\
\textsf{Author:} &\autherp
\end{tabular}
}
%\mline
{\footnotesize	
\tableofcontents }
}

%%%%%%%%%%%%%%%%%%%%%%%%%%%%%%%%%%
\def\pdf#1{{\rm PDF}_{#1}}
\def\tilco#1{\tilde{\co{#1}}}
\def\ndm#1{{\bf M}_{#1}(\{0,1\})}
\def\cdm#1{{\cal M}_{#1}(\{0,1\})}
\def\tndm#1{\tilde{{\bf M}}_{#1}(\{0,1\})}
%
%\def\openone%{\hbox{\upshape \small1\kern-3.3pt\normalsize1}}
%{\mathchoice
%{\hbox{\upshape \small1\kern-3.3pt\normalsize1}}
%{\hbox{\upshape \small1\kern-3.3pt\normalsize1}}
%{\hbox{\upshape \tiny1\kern-2.3pt\SMALL1}}
%{\hbox{\upshape \Tiny1\kern-2pt\tiny1}}}
%%%%%%%%%%%%% poor man's outline %%%%%%%%%%%%%%%%%%%
%\def\pmoutlinefnt#1{\setbox0=\hbox{#1}%
%   \setbox1=\hbox{\kern-.020em\copy0\kern-\wd0\kern.020em\copy0%
%   \kern-\wd0\kern.020em\copy0}$
%   \copy1\kern-\wd1\raise.020em\copy1\kern-\wd1\raise-.020em\copy1%
%   \color[rgb]{1,1,1}\kern-\wd0\kern-.020em\box0
%$}
%\def\openone{\mbox{\pmoutlinefnt{1}}}
\def\openone{\mbox{{\rm 1\hspace{-1mm}l}}}
\def\goh{{\goth h}}
%
%%%%%%%%% Cut from here %%%%%%%%%%
%\input comm.txt
%%%%%%%%% End of Cut %%%%%%%%%
%
%
\setcounter{section}{0}
\setcounter{footnote}{0}
\setcounter{page}{1}
\pagestyle{plain}

%
% Title
%
%%%%%%%%%%%%%%%%%%%%%%%%%%%%%%%%%%
\title{\titlep}
\author{\autherp\thanks{\emailp}
\\
\addressp}
\date{}
\maketitle
%%%%%%%%%%%%%%%%%%%%%%%%%
%
% Abstract
%
\begin{abstract}
Let ${\cal O}_n$ denote the Cuntz algebra for $2\leq n<\infty$.
With respect to a homogeneous embedding
of ${\cal O}_{n^m}$ into ${\cal O}_n$,
an extension of a Cuntz state on ${\cal O}_{n^m}$
to ${\cal O}_n$ is called a sub-Cuntz state,
which was introduced by Bratteli and Jorgensen.
We show (i) a necessary and sufficient condition of
the uniqueness of the extension,
(ii) the complete classification of pure sub-Cuntz states up to
unitary equivalence of their GNS representations, and 
(iii) the decomposition formula of a mixing sub-Cuntz state
into a convex hull of pure sub-Cuntz states.
Invariants of GNS representations of pure sub-Cuntz states
are realized as 
conjugacy classes of nonperiodic homogeneous unit vectors in a tensor-power vector space.
It is shown that this state parameterization satisfies
both the $U(n)$-covariance and the compatibility with a certain tensor product.
For proofs of main theorems,
matricizations of state parameters 
and properties of free semigroups are used.
\end{abstract}


%{\footnotesize	 \tableofcontents  }
%%%%%%%%%%%%%%%%%%%%
%
% Section 1
% 
\sftt{Introduction}
\label{section:first}
For a unital C$^*$-algebra $A$
and a unital C$^*$-subalgebra $B$ of $A$,
any state $\omega$ on $B$ has an 
{\it extension} $\tilde{\omega}$ on $A$,
that is,  $\tilde{\omega}$ is a state on $A$ which satisfies
$\tilde{\omega}|_B=\omega$
(\cite{Dixmier}, 2.10.1),
but it is not unique in general.
In this paper,
we completely classify extensions of a certain class of pure states 
on Cuntz algebras.
%with respect to certain embeddings between Cuntz algebras.
In consequence,
a new class of pure states on Cuntz algebras and 
the complete set of their invariants are given.
In this section, we show our motivation, definitions and main theorems.
Proofs will be given after $\S$ \ref{section:third}.

%%%%%%%%%%%%%%%%%%%%%%%%%%
%
% subsection 1.1
%
\ssft{Motivation}
\label{subsection:firstone}
In this subsection,  we make it clear that our aim of this study
against  a background of well-known representation theory,
and give a short survey. 
%%%%%%%%%%%%%%%%%%%%%%%%%%%%%
%
% subsubsection 1.1.1
% 
\sssft{Toward a representation theory of  C$^*$-algebras}
\label{subsubsection:firstoneone}
According to Kobayashi \cite{Kobayashi},
central problems of representation theory (except interactions with other
branches of mathematics)
%for a given group $G$ (or algebra $A$)
are listed as follows:
\begin{description}
\item[Problem 1]
Understanding irreducible representations.
Find and classify ``smallest" objects. 
There are the following subproblems:
\begin{itemize}
\item Construction of irreducible representations.
\item Finding a complete set of {\it invariants} of representations,
so that they can separate different irreducible representations from one another.
\item
Understanding these invariants.
\end{itemize}
\item[Problem 2]
Decompose a given representation into irreducible ones.
How is a given representation built from ``smallest" objects? 
\item[Problem 2-A]
Given an irreducible representation $\tau$ of a subgroup $G'$,
decompose the induced representation ${\rm Ind}_{G'}^G\tau$
into irreducibles of $G$.
\item[Problem 2-B]
Given an irreducible representation $\pi$ of $G$,
decompose the restriction $\pi|_{G'}$ into irreducibles of a subgroup $G'$.
The formula of the irreducible decomposition in this problem is called a
{\it branching law}
(e.g.,  the decomposition of tensor product representation).
\end{description}
Each problem is  more closely explained in the original text. 
Here a ``representation" means a representation of a group $G$.
%which may be assumed  of type I.
We wish to generalize 
Kobayashi's problems to the class of  algebras which includes  group algebras.

In general,
representations of C$^{*}$-algebras do not have unique decomposition
(up to
unitary equivalence) into sums or integrals of irreducibles \cite{Glimm}.
This is a difficulty to consider   
Problem 2 in the representation theory of C$^*$-algebras.
However,  it does not mean that 
every irreducible decomposition of a representation of a C$^*$-algebra
makes no sense.
If one chooses a good class of representations,
then Problem 2 can be treated satisfactorily.
For example, it is known that 
Cuntz algebras have such good classes of representations.

%%%%%%%%%%%%%%%%%%%%%%%%
%
% subsubsection 1.1.2
%
\sssft{States on Cuntz algebras}
\label{subsubsection:firstonetwo}
By Gel'fand-Naimark-Segal (=GNS) construction, 
the state theory of a C$^*$-algebra $A$ can be interpreted as  
the  (cyclic) representation theory of $A$ almost all.
Hence, we mainly consider (pure) states instead of (irreducible) 
representations in this paper.
For Cuntz algebras,
representations and states have been studied by many authors
\cite{ACE,BC,BJ1997,BJ,BJKW,BJO,BJP,Evans,FL,Gabriel,Jeong1999,
Jeong2005,GP0123,Laca1993,LS,Shin}
(see a specific survey in $\S$ 1 and $\S$ 2 of \cite{DHJ}),
but their classifications have not been finished yet.
The known most general approach was given in \cite{BJKW}.
The set of all states on a Cuntz algebra
is divided into two subsets, the set of finitely correlated states
and otherwise (= the set of infinitely correlated states) 
%\cite{BJ}
(see  $\S$ \ref{subsection:secondone}).
A finitely correlated state is characterized by
the existence of a finite-dimensional non-trivial $s_i^*$-invariant subspace 
of the GNS representation space
\cite{BJKW,DHJ}.
For example,
any Cuntz state (see $\S$ \ref{subsection:firsttwo}) is finitely correlated.
There exist both finitely and infinitely correlated 
vector states of permutative representations \cite{BJ}
(see $\S$ \ref{subsection:fourthtwo} and 
Example \ref{ex:permurep}).

We illustrate a rough classification of states on $\con$ 
($2\leq n<\infty$) as follows:

\def\state{
\put(-120,-100){
\begin{minipage}[t]{2in}
States on $\con$\\
($2\leq n<\infty$)
\end{minipage}
}
\put(210,-10){
\begin{minipage}[t]{2in}
infinitely correlated\\
states
\end{minipage}
}
\put(630,-5){\path(-30,0)(30,0)(0,0)(0,100)(30,100)}
\put(700,85){other}
\put(700,-10){quasi-free states}
\put(210,-310){
\begin{minipage}[t]{2in}
finitely correlated\\
states
\end{minipage}
}
\put(590,-300){\path(-30,0)(30,0)(0,0)(0,100)(30,100)}
\put(650,-310){{\bf 
\begin{minipage}[t]{2in}
sub-Cuntz\\
states
\end{minipage}
}}
\put(650,-210){other}
%\put(900,-110){permutative states}
\put(970,-310){Cuntz states}
\put(920,-300){\path(-30,0)(30,0)(0,0)(0,100)(30,100)}
\put(970,-210){other}
\put(955,-450){
\begin{minipage}[t]{2in}
Cuntz states\\
on $\co{n^m}$
\end{minipage}
}
\put(0,0){
\put(960,-550){\spline(0,100)(-150,150)(-150,180)}
\put(810,-359){\vector(0,1){0}}
\put(660,-450){{\sf extension}}
}
\put(170,-300){
\put(-30,200){\line(1,0){30}}
\path(30,0)(0,0)(0,300)(30,300)
}
}
\noindent
\thicklines
\setlength{\unitlength}{.1mm}
%\framebox{
\begin{picture}(1300,750)(-100,-550)
\put(-120,150){
\begin{minipage}[t]{2in}
%
% Figure 1.1
%
\begin{fig}
\label{fig:one}
\end{fig}
\end{minipage}
}
\put(0,0){\state}
\end{picture}

\noindent
About grounds in Figure \ref{fig:one},
see the proof of Fact \ref{fact:existence},
Lemma \ref{lem:omegasub}(i)
and Example \ref{ex:quasifree}.
Remark that Figure \ref{fig:one} is not true for the case of $\coni$
(see Proposition \ref{prop:infsub}).
Cuntz states are completely classified pure states with explicit complete invariants,
and are used to construct multiplicative isometries (\cite{TS08}, $\S$ 3)
and $R$-matrices (\cite{TS15}, $\S$ 3.2) (see also  \cite{PFO01,TS11}).
In this study, we select sub-Cuntz states as a target of complete classification
because they are natural generalizations of Cuntz states.
As well as Cuntz states,
it is expected that sub-Cuntz states have many applications.
Cuntz states and sub-Cuntz states will be explained explicitly  
in $\S$ \ref{subsection:firsttwo}.
Examples will be shown in $\S$ \ref{section:fourth}.
%%%%%%%%%%%%%%%%%%%
%
% subsubsection 1.1.3
%
\sssft{Branching laws of representations of Cuntz algebras}
\label{subsubsection:firstonethree}
We have mainly studied branching laws of representations  of Cuntz algebras
according to Kobayashi's Problem 2-B.
In \cite{PE01,PE02,PE03},
branching laws of permutative representations of Cuntz algebras
arising from endomorphisms  were computed
(see also \cite{Lawson2009}).
In \cite{AK03,AK02RR,IWF},
representations of fermions
were considered as restrictions of representations of $\co{2}$
by a certain embedding of the CAR algebra into $\co{2}$.
By using a certain set of embeddings between Cuntz algebras,
we defined a non-symmetric  tensor product of representations \cite{TS01}
(see $\S$ \ref{subsubsection:firstthreetwo}).
We showed the decomposition formula of this tensor product
of permutative representations
\cite{TS01,TS07}.
The set of all unitary equivalence classes 
of irreducible permutative representations of $\coni$
is one-to-one correspondence in the set of
all equivalence classes of  irrational numbers by modular transformations
\cite{CFR02}.
In this case,
finitely and infinitely correlated vector states 
associated with irreducible components are corresponded to
quadratic irrationals and otherwise, respectively.

Their common foundation is the representation
theory of Cuntz algebras.
Hence its development will be directly reflected in these subjects.

%%%%%%%%%%%%%%%%%%%%%%%%%%%%%%%%
%
% subsection 1.2
%
\ssft{Definition and main theorems}
\label{subsection:firsttwo}
In this subsection, we review the definition of sub-Cuntz state
by Bratteli-Jorgensen \cite{BJ},
and show our main theorems for $\con$ ($2\leq n<\infty$).
For the case of $\coni$, see Appendix \ref{section:appthree}.
Fix $2\leq n<\infty$.
Let $\con$ denote the {\it Cuntz algebra} with Cuntz generators
$s_1,\ldots,s_n$, that is,
$\con$ is a C$^{*}$-algebra which is universally generated by
$s_{1},\ldots,s_{n}$ which satisfy $s_{i}^{*}s_{j}=\delta_{ij}I$
for $i,j=1,\ldots,n$ and $\sum_{i=1}^{n}s_{i}s_{i}^{*}=I$ \cite{C}.

We review Cuntz state before sub-Cuntz state.
For any complex unit vector $z=(z_1,\ldots,z_n)\in {\Bbb C}^n$,
a state $\omega_z$ on $\con$ which satisfies
%
% Equation 1.1
%
\begin{equation}
\label{eqn:first}
\omega_z(s_j)=\overline{z_j}\quad\mbox{for all }j=1,\ldots,n,
\end{equation}
exists uniquely and is pure
where $\overline{z_i}$ denotes the complex conjugate of $z_i$.
The state $\omega_z$ is called the {\it Cuntz state} 
by $z$ \cite{BJ1997,BJ,BJP}.
GNS representations by $\omega_z$ and $\omega_y$
are unitarily equivalent if and only if $z=y$
(see Appendix \ref{section:apptwo}).
About equivalent definitions, see the case of $m=1$ in Theorem \ref{Thm:bj}.

For $m\geq 1$,
let ${\cal V}_{n,m}$ denote
the Hilbert space with an orthonormal basis
$\{e_J: J\in\{1,\ldots,n\}^m\}$,
that is, ${\cal V}_{n,m}=\ell^2(\{1,\ldots,n\}^m)$.
Let  $({\cal V}_{n,m})_1
:=\{z\in {\cal V}_{n,m}:\|z\|=1\}$.
%
% Definition 1.2
%
\begin{defi}
%\cite{BJ}
\label{defi:subcuntz}
For $z=\sum z_Je_J\in ({\cal V}_{n,m})_1$,
$\omega$ is a sub-Cuntz state on $\con$ by $z$ 
if $\omega$ is a state on $\con$ which satisfies the following equations:
%
% Equation 1.2
%
\begin{equation}
\label{eqn:subeqn}
\omega(s_{J})=\overline{z_{J}}
\quad \mbox{for all }J\in \{1,\ldots,n\}^m
\end{equation}
where $s_J:=s_{j_1}\cdots s_{j_m}$ when
$J=(j_1,\ldots,j_m)$,
and $\overline{z_J}$ denotes the complex conjugate of $z_J$.
In this case,
$\omega$ is called a sub-Cuntz state of order $m$.
\end{defi}

\noindent
This definition is equivalent to the original in \cite{BJ}
(see Theorem \ref{Thm:bj}).
%
% Fact 1.3
%
\begin{fact}(Existence)
\label{fact:existence}
For any $z\in ({\cal V}_{n,m})_1$,
a sub-Cuntz state by $z$ exists.
\end{fact}
%
% Proof 
%
\pr
Fix a bijection $f:\{1,\ldots,n^m\}\cong \{1,\ldots,n\}^m$.
Let $t_1,\ldots,t_{n^m}$ denote
the Cuntz generators of $\co{n^m}$.
Define the embedding $\hat{f}$ of $\co{n^m}$
into $\con$ by
$\hat{f}(t_i):=s_{j_1}\cdots s_{j_m}$
for $i\in\{1,\ldots,n^m\}$
when $f(i)=(j_1,\ldots,j_m)$.
We identify $\co{n^m}$ with $\hat{f}(\co{n^{m}})$ here. 
By definition,
$\omega$ is a sub-Cuntz state by $z$ if and only if 
$\omega$ is an extension of the Cuntz state $\omega_{\hat{z}}$
on $\co{n^m}$ to $\con$
where
$\hat{z}:=(z_{f(i)})_{i=1}^{n^m}\in ({\Bbb C}^{n^m})_1$.
Since an extension of $\omega_{\hat{z}}$ always exists,
the statement holds.
\qedh

\noindent
From the proof of Fact \ref{fact:existence},
a sub-Cuntz state is regarded as an extension of a Cuntz state.
Such an extension always exists but
it is not always unique.
We show a necessary and sufficient condition
of its uniqueness as follows.
Here we identify 
${\cal V}_{n,m}$ with $({\cal V}_{n,1})^{\otimes m}$
by the correspondence between bases
$e_{J}\mapsto e_{j_1}\otimes \cdots \otimes e_{j_m}$
for $J=(j_1,\ldots,j_m)\in\{1,\ldots,n\}^m$.
From this identification,
we obtain
${\cal V}_{n,m}\otimes {\cal V}_{n,l}={\cal V}_{n,m+l}$
 for any $m,l\geq 1$.
Then the following hold.
%
% Theorem 1.4
% 
\begin{Thm}
\label{Thm:maintwo}
Let $\omega$ be a sub-Cuntz state on $\con$ by $z\in ({\cal V}_{n,m})_1$.
\begin{enumerate}
%(i)
\item (Uniqueness)
$\omega$ is unique if and only if  $z$ is nonperiodic (or primitive \cite{Lothaire}),
that is, $z=x^{\otimes p}$ for some $x$ implies $p=1$.
In this case,
we write $\tilde{\omega}_z$ as $\omega$.
%(ii)
\item 
If $z$ is nonperiodic, then 
$\tilde{\omega}_z$ is pure.
%(iii)
\item
If $z$ is nonperiodic,
then 
%
% Equation 1.3
%
\begin{equation}
\label{eqn:omegasjsk}
\tilde{\omega}_z(s_Js_{K}^*)=\left\{
\begin{array}{ll}
\overline{z_{J}}\,z_{K}\quad &\mbox{when } |J|,|K|\in m{\Bbb Z}_{\geq 0},\\
\\
0 \quad& \mbox{when } |J|-|K|\not \in m{\Bbb Z},\\
\\
\overline{z_{J_1}}\,z_{K_1}
\disp{\sum_{|L|=m-|J_2|}\overline{z_{J_2L}}\,z_{K_2L}}
\quad &\mbox{otherwise}
\end{array}
\right.
\end{equation}
for $J,K\in\bigcup_{a\geq 1}\{1,\ldots,n\}^a\cup\{\emptyset\}$
where 
$|J|$ denotes the word length of $J$, 
$JK$ denotes the concatenation of $J$ and $K$,
$s_{\emptyset}:=I$,  $z_{\emptyset}:=1$ and
$z_{J}:=z_{J^{(1)}}\cdots z_{J^{(l)}}$ 
when $J=J^{(1)}\cdots J^{(l)}$
and $|J^{(i)}|=m$ for $i=1,\ldots,l$.
In the case of ``otherwise" in (\ref{eqn:omegasjsk}),
$J$ and $K$ satisfy
$J=J_1J_2$ and $K=K_1K_2$ such that
$|J_1|, |K_1| \in m{\Bbb Z}_{\geq 0}$ and 
$1\leq  |J_2|=|K_2|\leq m-1$.
\end{enumerate}
\end{Thm}

\noindent
From Theorem \ref{Thm:maintwo}(i),
if $z$ is periodic (= not nonperiodic), then $\omega$ is not unique.
In this case, all possibilities of sub-Cuntz states by $z$
are explicitly given as follows.
%
% Theorem 1.5
%
\begin{Thm}(Decomposition)
\label{Thm:periodic}
Let $p\geq 2$ and $z:=x^{\otimes p}$
for a nonperiodic element $x\in ({\cal V}_{n,m'})_1$.
If $\omega$ is a sub-Cuntz state on $\con$ by $z$,
then there exists $a=(a_1,\ldots,a_p)$ in 
$\Delta_{p-1}:=\{(b_1,\ldots,b_p)
\in {\Bbb R}^{p}:
\sum b_j=1,\,b_i\geq 0\mbox{ for all }i\}$
such that
$\omega$ has the following form:
%
% Equation 1.4
%
\begin{equation}
\label{eqn:omegatwo}
\omega=\sum_{j=1}^{p}a_j\omega_j
\end{equation}
where 
$\omega_j$ denotes the pure sub-Cuntz state by $e^{2\pi j\sqrt{-1}/p}x$.
%where $\zeta_p:=e^{2\pi \sqrt{-1}/p}$.
% denotes the $p$-th root of unity.
In (\ref{eqn:omegatwo}),
$(a_1,\ldots,a_p)$ is unique. 
\end{Thm}

\noindent
From Theorem \ref{Thm:periodic},
we see that a sub-Cuntz state by $z$ may be pure even if $z$ is periodic.
By combining Theorem \ref{Thm:maintwo}(ii) and 
Theorem \ref{Thm:periodic},
we obtain the following necessary and sufficient condition of 
the purity of sub-Cuntz state.
%
% Corollary 1.6
%
\begin{cor}(Purity)
\label{cor:pure}
For a sub-Cuntz state  $\omega$  by  $z\in ({\cal V}_{n,m})_1$,
$\omega$ is pure if and only if 
$\omega=\tilde{\omega}_x$
for some nonperiodic element $x\in ({\cal V}_{n,m'})_1$.
In this case,
$z=x^{\otimes p}$ for some $p\geq 1$.
\end{cor}
%
% Proof
%
\pr
($\Rightarrow$)
Assume that $\omega$ is pure.
If $z$ is nonperiodic,
then let $x:=z$.
If $z=v^{\otimes p}$ for some nonperiodic element $v$ and $p\geq 2$,
then $\omega=\sum_{j=1}^{p}a_j \omega_j$ 
from Theorem \ref{Thm:periodic}
where $\omega_j$ denotes the pure sub-Cuntz state by $e^{2\pi j\sqrt{-1}/p} v$.
Since $\omega$ is pure and $\omega_i\ne \omega_j$ when $i\ne j$,
there must exist $j$ such that $a_j=1$ and $\omega=\omega_j$.
Let $x:=e^{2\pi j\sqrt{-1}/p} v$.
Then  the statement holds.

\noindent
($\Leftarrow$) 
From Theorem \ref{Thm:maintwo}(ii), the statement holds.

From the above proofs,
we see that $z=x^{\otimes p}$ for $p\geq 1$.
\qedh

Next, we consider an equivalence of sub-Cuntz states.
%
% Theorem  1.7
%
\begin{Thm}(Equivalence)
\label{Thm:mainthree}
For $z,y\in \bigcup_{m\geq 1}({\cal V}_{n,m})_1$,
assume that both $z$ and $y$ are nonperiodic.
Then the following are equivalent:
\begin{enumerate}
%(i)
\item
GNS representations by $\tilde{\omega}_z$ and $\tilde{\omega}_y$
are unitarily equivalent. 
In this case,
we write $\tilde{\omega}_z\sim \tilde{\omega}_y$.
%(ii)
\item
\begin{enumerate}
%(a)
\item
$z=y$, or 
%(b)
\item
$z=x_1\otimes x_2$ and $y=x_2\otimes x_1$
for some $x_1, x_2\in \bigcup_{m\geq 1}({\cal V}_{n,m})_1$.
\end{enumerate}
In these cases, $z$ and $y$ are said to be conjugate (\cite{Lothaire}, $\S$ 1.3), 
and we write $z\sim y$.
\end{enumerate}
\end{Thm}

\noindent
Assume that 
both $z\in ({\cal V}_{n,m})_1$ and $y\in({\cal V}_{n,l})_1$
are nonperiodic.
If $m\ne l$, then $z\not\sim y$.
Hence $\tilde{\omega}_z\not\sim \tilde{\omega}_y$
from Theorem \ref{Thm:mainthree}.
From Theorem \ref{Thm:periodic},
two sub-Cuntz states of different orders
may be equivalent.
%
% Remark 1.8
%
\begin{rem}
\label{rem:first}
{\rm
\begin{enumerate}
%(i)
\item
By definition,
$\omega$ is a sub-Cuntz state of order $1$ 
if and only if $\omega$ is a Cuntz state.
From Theorem \ref{Thm:mainthree},
any Cuntz state is not equivalent to any sub-Cuntz state by 
a nonperiodic parameter $z\in ({\cal V}_{n,m})_1$ for $m\geq 2$.
%(ii)
\item
The restriction of any sub-Cuntz state on $\con$ by $z$ on
the UHF subalgebra $UHF_n:=C^*\{s_Js_K^*\in\con:|J|=|K|\}$ 
of $\con$ is always uniquely defined by $z$ 
from Theorem \ref{Thm:maintwo}(iii)
(see also \cite{BJ}, Proposition 5.1).
%Hence there is no problem of uniqueness on $UHF_n$.
%(iii)
\item
For $\omega$ in (\ref{eqn:omegatwo}),
we have the unique irreducible decomposition of the GNS representation $\pi$
by $\omega$ as follows:
%
% Equation 1.5
%
\begin{equation}
\label{eqn:pibig}
\pi=\bigoplus_{j=1}^{p}\hat{a}_j\pi_j,\quad
\hat{a}_{j}:=\left\{
\begin{array}{ll}
1\quad & (a_j> 0),\\
\\
0\quad & (a_j=0)
\end{array}
\right.
\end{equation}
where $\pi_j$ denotes the (irreducible) GNS representation
by $\omega_j$, and $\hat{a}_{j}$ means
the multiplicity coefficient of $\pi_j$ in $\pi$.
By Theorem \ref{Thm:mainthree},
$\pi_i\not\sim \pi_{j}$ when $i\ne j$
because $e^{2\pi i\sqrt{-1}/p}x\not\sim e^{2\pi j\sqrt{-1}/p}x$
when $i\ne j$.
Hence $\pi$ is multiplicity free.
In consequence,
the GNS representation by any sub-Cuntz state 
is multiplicity free,
and the class of GNS representations by sub-Cuntz states 
is closed with respect to the irreducible decomposition.
%(iv)
\item
We can verify that 
$\sim$ in Theorem \ref{Thm:mainthree}(ii)
is an equivalence relation.
Let $\approx$ denote the equivalence relation
on ${\cal V}_{n,m}=({\Bbb C}^n)^{\otimes m}$ 
by the action of the cyclic group ${\Bbb Z}/m{\Bbb Z}$
with respect to permutations of tensor components. 
Then $\sim$ does not coincide with $\approx$.
For example,
define three vectors in $({\cal V}_{n,3})_1$, $n\geq 3$ by
%
% Equation 1.6
%
\begin{equation}
\label{eqn:root}
\left\{
\begin{array}{rl}
z^{(1)}:=&\disp{e_1\otimes \frac{e_1\otimes e_2+e_2\otimes e_3}{\sqrt{2}},}\\
\\
z^{(2)}:=&\disp{\frac{e_1\otimes e_2+e_2\otimes e_3}{\sqrt{2}}\otimes e_1,}\\
\\
z^{(3)}:=&
\disp{\frac{e_2\otimes e_1\otimes e_1+e_3\otimes e_1\otimes e_2}{\sqrt{2}}.}
\end{array}
\right.
\end{equation}
Then  $z^{(1)}\approx z^{(2)}\approx z^{(3)}$,
but $z^{(1)}\sim z^{(2)}\not\sim z^{(3)}$.
\end{enumerate}
}
\end{rem}

%%%%%%%%%%%%%%%%%%%%%%%%%%%%%%%%
%
% subsection 1.3
%
\ssft{Naturalities of state parameterization}
\label{subsection:firstthree}
In Theorem \ref{Thm:maintwo},
we introduced a parametrization of pure sub-Cuntz states:
%
% Equation 1.7
%
\begin{equation}
\label{eqn:parameter}
z\longmapsto \tilde{\omega}_z.
\end{equation}
%
%for a nonperiodic parameter $z\in ({\cal V}_{n,m})_1$.
In this subsection,
we show how natural this parameterization is.
%%%%%%%%%%%%%%%%%%%%%%%%%%%%%
%
% subsection 1.3.1
%
\sssft{$U(n)$-covariance}
\label{subsubsection:firstthreeone}
For convenience,
%In order to treat only the unique extension case, 
we introduce some symbols as follows.
%
% Corollary 1.9
%
\begin{cor}
\label{cor:sc}
Define 
\[\begin{array}{rl}
{\goth N}_{n}:=&
\{z\in \bigcup_{m\geq 1}({\cal V}_{n,m})_1:z\mbox{ is nonperiodic}\},\\
\\
{\goth I}_{n}:=&
\{z\in \bigcup_{m\geq 1}({\cal V}_{n,m})_1:z\mbox{ is indecomposable}\},\\
\\
{\cal P}_n:&\mbox{the set of all pure states on $\con$},\\
\\
{\cal P}_{n,sub}:&\mbox{the set of all pure sub-Cuntz states on $\con$},\\
\\
{\rm Spec}\con:&\mbox{
\begin{minipage}[t]{4in}
the set of all unitary equivalence classes\\
of irreducible representations of $\con$
\end{minipage}
}\\
\end{array}
\]
where $z$ is said to be {\it indecomposable} if $z$ can not be 
written as $z_1\otimes z_2$ for any $z_1,z_2$.
% When $m=2$, a decomposable $z$ is said to be {\it pure or simple}.
% wikipedia
Then the following hold:
\begin{enumerate}
%(i)
\item
The map $q:{\goth N}_n\to {\cal P}_{n,sub};\, q(z):=\tilde{\omega}_z$,
is bijective.
%(i)
\item
The map $r:{\goth I}_n\to{\rm Spec}\con;\,r(z):=[\pi_z]$,
is injective
where $[\pi_z]$ denotes the unitary equivalence class
of the GNS representation $\pi_z$ by $\tilde{\omega}_z$.
\end{enumerate}
\end{cor}
%
% Proof
%
\pr
(i)
From Corollary \ref{cor:pure},
the statement holds.

\noindent
(ii)
From Theorem \ref{Thm:mainthree},
if both $z$  and $y$ are indecomposable,
then $\tilde{\omega}_z\sim \tilde{\omega}_y$  if and only if $z=y$.
From this, the statement holds.
\qedh

\noindent
From Theorem \ref{Thm:mainthree},
${\goth N}_n/\!\!\!\sim\,\,\,\cong \,\,{\cal P}_{n,sub}/\!\!\!\sim$\,
and ${\goth I}_n/\!\!\!\sim \,\,={\goth I}_n$.
The parameter set ${\goth I}_n$
can be regarded as the set 
%$\{f\in {\Bbb C}[x_1,\ldots,x_n]_{hom, irr}:\|f\|=1\}$
of non-commutative homogeneous 
irreducible polynomials in $n$-variables with the norm $1$ \cite{Lang}.
%where $\|f\|:=\{\sum_{J}|f_J|^2\}^{1/2}$
%when $f=\sum_{J}f_Jx_J$.

We show a naturality of the parameterization in Corollary \ref{cor:sc}
with respect to the standard unitary group action $\alpha$ on $\con$,
which is defined as
%
% Equation 1.8
%
\begin{equation}
\label{eqn:alphagsi}
\alpha_g(s_i):=
\sum_{j=1}^{n}g_{ji}s_{j}
\quad(i=1,\ldots,n,\,g=(g_{ij})\in U(n)).
\end{equation}
Define the dual action $\alpha^*$ of $\alpha$
on the dual $\con^*$ of $\con$ by
$\alpha^*_g(f):=f\circ \alpha_{g^*}$ for $f\in \con^*$ and $g\in U(n)$.
Especially, $\alpha^*_g({\cal P}_n)={\cal P}_n$ for all $g$.

Let  $\gamma$ denote  the standard action of $U(n)$ on
${\cal V}_{n,1}={\Bbb C}^n$, that is,
%
% Equation 1.9
%
\begin{equation}
\label{eqn:deeltage}
\gamma_ge_{i}:=\sum_{j=1}^{n}g_{ji}e_j\quad(i=1,\ldots,n,\,g\in U(n)).
\end{equation}
%
%Let $\gamma^{\otimes m}$ denote the $m$-times tensor power
%of $\gamma$, that is,
%$\gamma^{\otimes m}_g:=(\gamma_g)^{\otimes m}$.
Since $\gamma^{\otimes m}_g:=(\gamma_g)^{\otimes m}$ is a unitary,
$\gamma^{\otimes m}_g(({\cal V}_{n,m})_1)=({\cal V}_{n,m})_1$
for all $g\in U(n)$.
Define the action $\Gamma$ of $U(n)$ 
on $\bigcup_{m\geq 1}({\cal V}_{n,m})_1$ by
$\Gamma_gz:=\gamma^{\otimes m}_g z$ when $z\in ({\cal V}_{n,m})_1$.
Remark that if $z$ is nonperiodic ({\it resp.} indecomposable),
then $\Gamma_gz$ is also nonperiodic
({\it resp.} indecomposable)  for any $g\in U(n)$.
%
% Proposition  1.10
%
\begin{prop}($U(n)$-covariance)
\label{prop:cova}
Let ${\goth N}_n$ be as in Corollary \ref{cor:sc}.
For any $g\in U(n)$ and $z\in {\goth N}_n$,
%
% Equation 1.11
%
\begin{equation}
\label{eqn:alphag}
\alpha^*_g(\tilde{\omega}_{z})
=\tilde{\omega}_{\Gamma_gz}.
\end{equation}
That is, the parameterization 
$z\mapsto \tilde{\omega}_z$ is $U(n)$-covariant.
\end{prop}
%
% Proof
%
\pr
Assume $z\in ({\cal V}_{n,m})_1\cap {\goth N}_n$.
By definition,
we can verify 
$\alpha^*_g(\tilde{\omega}_{z})(s_J)
=\overline{(\gamma^{\otimes m}_gz)_J}
=\overline{(\Gamma_gz)_J}$
for all $J\in \{1,\ldots,n\}^m$ and $g\in U(n)$.
Since $\Gamma_gz$ is nonperiodic,
$\alpha^*_g(\tilde{\omega}_{z})$
coincides with 
$\tilde{\omega}_{\Gamma_gz}$
from Theorem \ref{Thm:maintwo}(i).
\qedh

\noindent
In other words,
%$q\circ \Gamma_g=\alpha_g^*\circ q$ for any $g\in U(n)$
%for $q$ in Corollary \ref{cor:sc}(i).
%Geometrically, 
%$q$ preserves any $U(n)$-orbit from ${\goth N}_n$
%to ${\cal P}_n$.
$q$ in Corollary \ref{cor:sc}(i) is an isomorphism between
two dynamical systems
$({\goth N}_n,\Gamma,U(n))$ and $({\cal P}_{n,sub},\alpha^*,U(n))$.
%are isomorphic as dynamical systems.

%%%%%%%%%%%%%%%%%%%%%%%%%%%%%
%
% subsection 1.3.2
%
\sssft{Compatibility with $\varphi$-tensor product}
\label{subsubsection:firstthreetwo}
%%%%%%%%%%%%%%%%%%%%
%
In \cite{TS01},
we introduced a non-symmetric tensor product of 
states on Cuntz algebras.
In this subsection,
we show tensor product formulas of sub-Cuntz states.
%From this,
%one more naturality of the parametrization of sub-Cuntz states
%in Corollary \ref{cor:sc} is shown.

We review definitions in \cite{TS01}.
Let $s_1^{(n)},\ldots,s_{n}^{(n)}$ denote Cuntz generators of $\con$.
For $2\leq n,n'<\infty$,
define the unital 
$*$-embedding $\varphi_{n,n'}$ of $\co{nn'}$ into $\con\otimes \co{n'}$ by
$\varphi_{n,n'}(s_{n'(i-1)+j}^{(nn')}):= s_{i}^{(n)}
\otimes s_{j}^{(n')}$ for $i\edot,\,j=1,\ldots,n'$.
Let ${\cal S}_n$ denote the set of 
all states on $\con$.
For $(\omega_1,\omega_2)\in {\cal S}_n\times {\cal S}_{n'}$,
the {\it $\varphi$-tensor product} 
$\omega_{1}\otimes_{\varphi} \omega_{2}\in{\cal S}_{nn'}$ is  defined by 
%
% Equation 1.12
%
\begin{equation}
\label{eqn:definition}
\omega_{1}\ptimes\omega_{2}:=(\omega_{1}\otimes \omega_{2})
\circ \varphi_{n,n'}.
\end{equation}

\noindent
Then $\ptimes$ is associative.
Hence the set $\bigcup_{n\geq 2}{\cal S}_n$ is a semigroup with the product
$\ptimes$.
Furthermore, the following holds.
%
% Proposition 1.11
%
\begin{prop}
\label{prop:quasitensor}
The $\varphi$-tensor product  of any two sub-Cuntz states
is also a sub-Cuntz state,
that is,
the set of all sub-Cuntz states is closed with respect to 
$\ptimes$.
\end{prop}

For $J=(j_1,\ldots,j_m)\in \{1,\ldots,n\}^{m}$
and 
$K=(k_1,\ldots,k_m)\in \{1,\ldots,n'\}^{m}$,
define
$J\boxtimes K=(l_1,\ldots,l_m)\in\{1,\ldots,nn'\}^m$ by
$l_t:=n'(j_t-1)+k_t$ for $t=1,\ldots,m$.
For $z=\sum z_Je_J\in {\cal V}_{n,m}$
and $y=\sum y_Ke_K\in {\cal V}_{n',m}$,
define $z\boxtimes y\in {\cal V}_{nn',m}$ by
%
% Equation 1.12
%
\begin{equation}
\label{eqn:boxtimes}
z\boxtimes y:=\sum_{L\in\{1,\ldots,nn'\}^m}(z\boxtimes y)_L\,e_{L},\quad 
(z\boxtimes y)_L:=z_{J}y_{K}
\end{equation}
where 
$J\in\{1,\ldots,n\}^m$
and $K\in\{1,\ldots,n'\}^m$
are uniquely defined as  $J\boxtimes K=L$.
By definition,
$\|z\boxtimes y\|=\|z\|\cdot \|y\|$.
If $\|z\|=\|y\|=1$,
then $\|z\boxtimes y\|=1$.
Remark that 
for any $z,y\in {\cal V}_{n,m}$,
$z\otimes y\in {\cal V}_{n,2m}$
and 
$z\boxtimes y\in {\cal V}_{n^2,m}$.
Clearly, ${\cal V}_{n,2m}\cong {\cal V}_{n^2,m}$,
but we distinguish $\otimes$ from $\boxtimes$ here.

In addition to $\boxtimes$,
we define a new operation.
For $z=\sum z_Je_J\in {\cal V}_{n,m}$
and $y=\sum y_Ke_K\in {\cal V}_{n',l}$,
define
%
% Equation 1.13
%
\begin{equation}
\label{eqn:albe}
z*y:=z^{\otimes \alpha}\boxtimes y^{\otimes \beta}\in
{\cal V}_{nn',d}
\end{equation}
where $d,\alpha,\beta\geq 1$ are uniquely chosen such that
$d:=\alpha m=\beta l$ is the least common multiple of $m$ and $l$.
Especially,
if $m=l$, then $\alpha=\beta=1$, $d=m$ and $z*y=z\boxtimes y$.
If $\|z\|=\|y\|=1$, then $\|z* y\|=1$.
About examples of these operations,
see \cite{TS01}.
%
% Proposition  1.12
%
\begin{prop}
\label{prop:tensor}
Assume that both
$z\in ({\cal V}_{n,m})_1$ and $y\in ({\cal V}_{n',l})_1$
are nonperiodic.
Then the following hold:
\begin{enumerate}
%(i)
\item
$z*y$ is nonperiodic.
%(ii)
\item (Tensor product formula)
$\tilde{\omega}_z\ptimes \tilde{\omega}_y=\tilde{\omega}_{z*y}$.
Especially,
if $m=l$,
then 
$\tilde{\omega}_z\ptimes \tilde{\omega}_y=\tilde{\omega}_{z\boxtimes y}$.
\end{enumerate}
\end{prop}

\noindent
Let $\omega_z$ denote the Cuntz state on $\con$ by $z\in ({\Bbb C}^n)_1$.
As special cases of 
Proposition \ref{prop:cova} and 
Proposition \ref{prop:tensor}(ii),
the following hold:
%
% Equation 1.14
%
\begin{equation}
\label{eqn:propro}
\alpha_g^*(\omega_z)=\omega_{gz},\quad
\omega_z\ptimes \omega_y=\omega_{z\boxtimes y}
\end{equation}
for  any $z\in ({\Bbb C}^n)_1$,
$y\in ({\Bbb C}^{n'})_1$ and $g\in  U(n)$
where $gz=\gamma_g z$. 
%Especially, the set of all Cuntz states
%is closed with respect to both the $U(n)$-action and $\ptimes$.

%\noindent
%Remark that Proposition \ref{prop:tensor} does not mean
%tensor product of GNS representations by $\tilde{\omega}_z$
%and $\tilde{\omega}_y$ are unitarily equivalent to
%that of $\tilde{\omega}_{z* y}$.

%\noindent
From 
Proposition \ref{prop:tensor}(ii), the operation $*$ is associative
because $\ptimes$ is associative.
Proposition \ref{prop:tensor}(i) means that 
${\goth N}_*:=\bigcup_{n\geq 2}{\goth N}_n$ is
a semigroup with the product $*$.
Furthermore, the following holds
from Corollary \ref{cor:pure}.
%
% Corollary 1.13
%
\begin{cor}
\label{cor:natural}
For $q$ and ${\cal P}_{n,sub}$ in Corollary \ref{cor:sc},
the set ${\cal P}_{*,sub}:=\bigcup_{n\geq 2}{\cal P}_{n,sub}$
is a semigroup with the product $\ptimes$,
and $q$ can be extended to an isomorphism 
between $({\goth N}_*,*)$ onto $({\cal P}_{*,sub},\ptimes)$.
\end{cor}

\noindent
Corollary \ref{cor:natural} means the second naturality of 
the state parametrization $q$.
%For Kobayashi's Problem 1 ``Understanding these invariants",
%the $U(n)$-covariance  and  
%the compatibility with $\varphi$-tensor product 
%of the state parametrization are our answers.

The paper is organized as follows:
In $\S$ \ref{section:second},
we will review known results and prepare tools to prove main theorems.
In $\S$ \ref{subsection:secondtwo},
a matricization of state parameter will be introduced.
In $\S$ \ref{section:third},
we will prove main theorems.
In $\S$ \ref{section:fourth},
we will show examples.
In $\S$ \ref{subsection:fourthone},
we will show sub-Cuntz states of order $2$.
In $\S$ \ref{subsection:fourthtwo},
sub-Cuntz states associated with 
permutative representations will be explained.
In $\S$ \ref{subsection:fourththree},
examples of non-sub-Cuntz states will be shown.

%%%%%%%%%%%%%%%%%%%%%%%%%%%%
%
% Section 2
%
\sftt{Preparations}
\label{section:second}
%
%In this section, we review known results and prepare tools to prove main theorems.
%%%%%%%%%%%%%%%%%%%
%
% subsection  2.1
%
\ssft{Finitely correlated states on $\con$}
\label{subsection:secondone}
We start from general properties of extensions of states.
%
% Proposition  2.1
%
\begin{prop}
\label{prop:stateext}
(\cite{Dixmier}, 2.10.1)
%p58, 
%
% K-R I, 4.3.13
Let $A$ be a C$^*$-algebra with unit $I$, and $B$ a C$^*$-subalgebra
of $A$ such that $I\in B$.
Then the following hold:
\begin{enumerate}
%(i)
\item
Every state on $B$ can be extended to a state on $A$.
%(ii)
\item
Every pure state on $B$ can be extended to a pure state on $A$.
Especially,
if its extension is unique, then it is pure.
\end{enumerate}
\end{prop}

\noindent
The existence of sub-Cuntz state is assured by 
Proposition \ref{prop:stateext}(i).
If it is  unique, then 
its purity  is assured by Proposition \ref{prop:stateext}(ii).

%
% Definition 2.2
%
\begin{defi}
\label{defi:finitely}
(\cite{BJ1997,BJKW})
A state $\omega$ on $\con$ ($2\leq n\leq \infty$) is 
said to be finitely correlated
if $\dim {\cal K}<\infty$
where 
${\cal K}:=\overline{{\rm Lin}\langle \{\pi(s_J)^*\Omega\in {\cal H}:J\}\rangle }$
and $({\cal H},\pi,\Omega)$ denotes the GNS representation by $\omega$.
If not, $\omega$ is said to be infinitely correlated. 
\end{defi}
%
%Both finitely  and infinitely correlated states appear as vector states of 
%irreducible components in  the shift representation of $\con$ \cite{BJ}
%(see also $\S$ \ref{subsection:fourthtwo} and $\S$ \ref{subsection:fourththree}).

Next, we show equivalent definitions of sub-Cuntz state as follows.
%
% Theorem 2.3
% 
\begin{Thm}
\label{Thm:bj}
Fix $m\geq 1$.
Let $\omega$ be a state on $\con$ with the GNS representation
$ ({\cal H}, \pi,\Omega)$.
For $z=\sum z_Je_J\in ({\cal V}_{n,m})_1$,
the following conditions are equivalent:
\begin{enumerate}
%(i)
\item
$\omega$ is a sub-Cuntz state by $z$.
%(ii)
\item
$\Omega=\pi(s(z))\Omega$
where $s(z):=\sum z_Js_J$.
%(iii)
\item
$\pi(s_{J})^*\Omega=z_{J}\Omega$
for all $J\in \{1,\ldots,n\}^m$.
\end{enumerate}
\end{Thm}
%
% Proof
%
\pr
From Proposition 5.1 of \cite{BJ},
(ii) and (iii) are equivalent.
By the definition of $({\cal H},\pi,\Omega)$,  (iii) implies (i). 
From (i), we have $\sum z_{J}\omega(s_{J})=1$.
This implies 
$\|\Omega-\sum z_{J}\pi(s_{J})\|=0$.
Hence (ii) holds.
\qedh

\noindent
In Definition 5.7 of \cite{DHJ}, 
a cyclic representation of $\con$ which satisfies 
equations in Theorem \ref{Thm:bj}(iii) with $m=1$
is called a {\it generic representation}.
%
% Lemma 2.4
%
\begin{lem}
\label{lem:omegasub}
\begin{enumerate}
%(i)
\item
When $n<\infty$,
any sub-Cuntz state on $\con$ is finitely correlated.
%(ii)
\item
If $\omega$ is a sub-Cuntz state with 
the GNS representation $({\cal H},\pi,\Omega)$,
then 
${\cal H}=\overline{{\rm Lin}\langle\{\pi(s_J)\Omega:J\}\rangle}$.
\end{enumerate}
\end{lem}
%
% Proof
%
\pr
Assume that $\omega$ is a sub-Cuntz state on $\con$
by $z=\sum z_Je_J\in ({\cal V}_{n,m})_1$.

\noindent
(i)
By Theorem \ref{Thm:bj}(iii),
$\dim {\rm Lin}\langle
\{\pi(s_J^*)\Omega:J\}
\rangle\leq 
\sum_{l=0}^{m-1}\dim {\rm Lin}\langle\{\pi(s_J^*)\Omega:|J|=l\}
\leq \sum_{l=0}^{m-1}n^{l}<\infty$.

\noindent
(ii)
Since $\con$ is spanned by $\{s_{J}s_K^*:J,K\}$,
${\cal H}$ is spanned by $\{\pi(s_{J}s_K^*)\Omega:J,K\}$.
From 
Theorem \ref{Thm:bj}(ii),
$\pi(s_{J}s_K^*)\Omega
=\pi(\,s_{J}s_K^*(s(z))^{l}\,)\Omega$
for any $l\geq 1$.
Therefore 
$\pi(s_{J}s_K^*)\Omega\in 
{\rm Lin}\langle\{\pi(s_{J'})\Omega:J'\}\rangle$
for any $J,K$.
This implies the statement.
\qedh

\noindent
Remark that Lemma \ref{lem:omegasub}(i) does not hold
for $\coni$ (see Proposition \ref{prop:infsub}).

%%%%%%%%%%%%%%%%%%%%%%%%%%%%
%
% subsection 2.2
% 
\ssft{Matricization of state parameter}
\label{subsection:secondtwo}
Assume $m\geq 2$.
In this subsection, we introduce operators
associated with an element $z\in {\cal V}_{n,m}$.
%In other words,
%it is a matrix representation of state parameter.

For $x\in {\cal V}_{n,a}$ 
and $y\in {\cal V}_{n,b}$ with $a,b\geq 1$,
define the operator $x\otimes y^*$ from ${\cal V}_{n,b}$ to ${\cal V}_{n,a}$ 
by $(x\otimes y^*)v:=\langle y|v\rangle x$
for $v\in {\cal V}_{n,b}$.
We generalize this as follows.
For $z=\sum z_Me_M\in {\cal V}_{n,m}$ and $1\leq a\leq m-1$,
define the operator
$T_{a}(z)$ from ${\cal  V}_{m,a}$ to ${\cal V}_{n,m-a}$ by
%
% Equation 2.1
%
\begin{equation}
\label{eqn:ta}
T_{a}(z)e_K:=\sum_{|J|=m-a}z_{JK}e_{J}\quad(K\in\{1,\ldots,n\}^{a}).
\end{equation}
In other words, 
$T_{a}(z)v=\sum_{|J|=m-a}\langle \,\bar{z}\,|e_{J}\otimes v\rangle e_J$
for $v\in {\cal V}_{n,a}$, or 
$T_a(z)=\sum_{J,K}z_{JK}\,e_J\otimes e_K^*$
where $\overline{z}:=\sum \overline{z_M}\,e_M$.
The operator $T_a(z)$ is called the {\it matricizing (matricization)} \cite{ES},
{\it unfolding} \cite{LMV} or {\it flattening} \cite{VT} of a tensor $z\in
{\cal V}_{n,m}=({\Bbb C}^n)^{\otimes m}$.
Especially,
$T_a(x\otimes y)=x\otimes \overline{y}^*$
for any $x\in {\cal V}_{n,m-a}$ and $y\in {\cal V}_{n,a}$.

In the case of $m=2$, 
$T_1(z)$ is identified with the matrix representation
 $(z_{ij})\in M_n({\Bbb C})$
of a $2$-tensor $z=\sum_{i,j=1}^{n}z_{ij}e_{ij}\in {\cal V}_{n,2}$ by definition.
In general case,
by the identification ${\rm Hom}_{{\Bbb C}}({\cal V}_{n,a},{\cal V}_{n,m-a})$
with the set $M_{n^{m-a},n^{a}}({\Bbb C})$
of all $n^{m-a}\times n^a$ matrices,
$T_a$ is regarded as the following  mapping:
%
% Equation 2.2
%
\begin{equation}
\label{eqn:matrix}
({\Bbb C}^{n})^{\otimes m}={\cal V}_{n,m}\ni z\longmapsto 
T_a(z)\in M_{n^{m-a},n^{a}}({\Bbb C}).
\end{equation}

In order to show properties of $T_a(z)$,
we review operator norms as follows.
Let ${\cal H}$ and ${\cal K}$ be Hilbert spaces.
For a  bounded linear operator $A$  from ${\cal H}$ to ${\cal K}$,
the {\it uniform norm} and the {\it Hilbert-Schmidt norm} of $A$ 
are defined as 
$\|A\|:=\sup_{x\in{\cal H},\,\|x\|=1}\|Ax\|$ and 
$\|A\|_2:=({\rm tr}A^*A)^{1/2}$, respectively \cite{Bhatia,DS2,GV}.
Then $\|A\|\leq \|A\|_2$.
Furthermore,
$\|A\|=\|A\|_2\ne 0$ if and only if  
there exist $y\in {\cal H}$ and $x\in {\cal K}$
such that $x,y\ne 0$ and $A=x\otimes y^*$.
%
% Fact 2.5
%
\begin{fact}
\label{fact:operator}
Let $z\in {\cal V}_{n,m}$.
\begin{enumerate}
%(i)
\item
For $1\leq a\leq m-1$,
$\|T_{a}(z)\|\leq \|z\|$.
%(ii)
\item
If $z\ne 0$, then  
$\|T_{a}(z)\|=\|z\|$
if and only if 
$z=x\otimes y$ for some $x\in {\cal V}_{n,m-a}$
and $y\in {\cal V}_{n,a}$.
\end{enumerate}
\end{fact}
%
% Proof
%
\pr
(i)
From the inequality of norms and (\ref{eqn:ta}),
%
% Equation 2.3
%
\begin{equation}
\label{eqn:ineqt}
\|T_{a}(z)\|\leq \|T_{a}(z)\|_2
=\Bigl\{\sum_{|K|=a}\,\sum_{|J|=m-a}|z_{JK}|^2\Bigr\}^{1/2}
=\Bigl\{\sum_{|M|=m}|z_{M}|^2\Bigr\}^{1/2}
=\|z\|.
\end{equation}

\noindent
(ii)
Assume $\|T_a(z)\|=\|z\|$.
From the proof of (i),  $\|T_{a}(z)\|=\|z\|=\|T_{a}(z)\|_2$.
Hence $T_a(z)=x\otimes w^*$
for some $x\in {\cal V}_{n,m-a}$
and $w\in {\cal V}_{n,a}$.
This implies
$z_{JK}=x_J\overline{w_K}
=(x\otimes \overline{w})_{JK}$
for $J\in \{1,\ldots,n\}^{m-a}$
and $K\in \{1,\ldots,n\}^{a}$.
By taking $y:=\overline{w}$,
$z=x\otimes y$.
The inverse direction holds from (\ref{eqn:ta}).
\qedh

The following is one of key lemmas to prove main theorems.
%
% Lemma 2.6
%
\begin{lem}
\label{lem:newone}
Let $m\geq 2$.
\begin{enumerate}
%(i)
\item
For $X,Y\in ({\cal V}_{n,m})_1$,
if there exists 
a nonzero vector $v\in {\cal V}_{n,a}$ 
for some $1\leq a\leq m-1$
which satisfies
%
% Equation 2.6
%
\begin{equation}
\label{eqn:uazetad}
v=c\, T_{m-a}(\overline{X})T_{a}(Y)v
\end{equation}
for some $c\in U(1):=\{c'\in {\Bbb C}:|c'|=1\}$,
then there exist $x_1\in ({\cal V}_{n,a})_1$ 
and $x_2\in ({\cal V}_{n,m-a})_1$
such that $X=x_1\otimes  x_2$ and $Y=\overline{c}\,x_2\otimes x_1$.
%(ii)
\item
In addition to (i), if $X=Y$,
then $c=1$ and 
%$X=x_1\otimes  x_2= x_2\otimes x_1$.
$X$ is periodic.
\end{enumerate}
\end{lem}
%
% Proof
%
\pr
(i)
From Fact \ref{fact:operator}(i),
$\|T_b(z)\|\leq \|z\|=1$ 
for $z=X,Y$ and $b=a,m-a$.
Since $v\ne 0$, we obtain
$\|T_a(Y)\|=\|T_{m-a}(X)\|=\|T_{m-a}(\overline{X})\|=1
=\|X\|=\|Y\|$
from (\ref{eqn:uazetad}).
From these and Fact \ref{fact:operator}(ii),
%
% Equation 2.7
%
\begin{equation}
\label{eqn:exchange}
X=x_1\otimes x_2, \quad Y=x_1'\otimes x_2'
\end{equation}
for some 
$x_1,x_2'\in ({\cal V}_{n,a})_1$
and
$x_2,x_1'\in ({\cal V}_{n,m-a})_1$.
By substituting (\ref{eqn:exchange}) into (\ref{eqn:uazetad}),
%
% Equation 2.7
%
\begin{equation}
\label{eqn:uazetac}
v=c\,\overline{x_1}\langle x_2|x_1'\rangle
\langle \overline{x_2'}| v\rangle.
\end{equation}
From (\ref{eqn:uazetac}),
$\langle \overline{x'_2}|v\rangle =c\langle \overline{x'_2}| \overline{x_1}\rangle
\,
\langle x_2|x_1'\rangle
\langle \overline{x_2'}| v\rangle$.
By $v\ne 0$ and (\ref{eqn:uazetac}),
$\langle \overline{x_2'}| v\rangle\ne 0$.
Hence 
$1 =c\langle \overline{x'_2}| \overline{x_1}\rangle
\,
\langle x_2|x_1'\rangle
=\langle \overline{c}x_2\otimes x_1 | x_1'\otimes x_2'\rangle$.
From Lemma \ref{lem:trivial},
$x_1'\otimes x_2'=\overline{c}x_2\otimes x_1$.
From this and  (\ref{eqn:exchange}),
the statement holds.

\noindent
(ii)
From (i),
$\overline{c}x_2\otimes x_1=Y=X=x_1\otimes x_2$.
Applying 
Corollary \ref{cor:projective}(i) to this, 
the statement holds.
\qedh

%%%%%%%%%%%%%%%%%%%%%%%
%
% subsection 2.3
%
\ssft{Reduction of problems}
\label{subsection:secondthree}
For a given sub-Cuntz state $\omega$ of order $m\geq 2$,
Definition \ref{defi:subcuntz}
determines only special values $\{\omega(s_{J}):|J|=m\}$,
but not all values $\{\omega(s_{J}s_{K}^*):J,K\}$.
From Theorem \ref{Thm:bj},
%By using these relations, 
we can reduce the uniqueness problem of $\omega$ to 
the problem to determine the smaller set $\{\omega(s_{J}):1\leq |J|\leq m-1\}$.
%
% Lemma 2.7
%
\begin{lem}
\label{lem:redthree}
For $z\in ({\cal V}_{n,m})_1$ and $y\in ({\cal V}_{n,l})_1$,
we introduce the following assumptions:\\
{\bf Assumption E}:
``$\con$ acts on a Hilbert space ${\cal H}$ with
two unit vectors $\Omega_z$ and $\Omega_y$  which satisfy
%
% Equation 2.8
%
\begin{equation}
\label{eqn:definitioneq}
s(z)\Omega_z=\Omega_z,\quad s(y)\Omega_y=\Omega_y
\end{equation}
where $s(z):=\sum z_Js_J\in\con$ for $z=\sum z_Je_J$."\\
{\bf Assumption EC}:
Assumption E with the cyclicity of both $\Omega_z$ and $\Omega_y$.

Assume Assumption E for $z$ and $y$.
Define the linear functional $\omega_{z,y}$ on $\con$ by
%
% Equation  2.9
%
\begin{equation}
\label{eqn:omegazyb}
\omega_{z,y}:=\langle \Omega_z|(\cdot)\Omega_y\rangle.
\end{equation}
Then the following hold.
\begin{enumerate}
\item
For $J,K\in \bigcup_{k\geq 0}\{1,\ldots,n\}^k$,
assume that
$J=J_1J_2, K=K_1K_2$,
$|J_1|=ma$, $|K_1|=lb$ for some
$a,b\geq 0$, 
and $|J_2|=i,|K_2|=i'$, 
$0\leq i\leq m-1$ and $0\leq i'\leq l-1$.
\begin{enumerate}
%(a)
\item
If $(i,i')=(0,0)$,
then 
$\omega_{z,y}(s_{J}s_K^*)=\overline{z_J}\,y_K\,\omega_{z,y}(I)$.
%(b)
\item
If $(i,i')\ne (0,0)$
and $\alpha m-i=\beta l-i'$ for some $\alpha,\beta\geq 1$,
then
%
% Equation 2.10
%
\begin{equation}
\label{eqn:sjskstar}
\omega_{z,y}(s_{J}s_{K}^*)
=
\overline{z_{J_1}}\,y_{K_1}
\omega_{z,y}(I)\,
\sum_{|A|=\alpha m-i}
\overline{(z^{\otimes \alpha})_{J_2 A}}\,
(y^{\otimes \beta})_{K_2 A}.
\end{equation}
%
%(c)
\item
If $(i,i')\ne (0,0)$ 
and $\alpha m-i\ne \beta l-i'$ for any $\alpha,\beta\geq 1$,
then 
$\omega_{z,y}(s_{J}s_{K}^*)\in {\cal W}_{z,y}
:={\rm Lin}\langle
\{\omega_{z,y}(s_{M}),
\omega_{z,y}(s_{L}^*):
1\leq |M|\leq m-1,
1\leq |L|\leq l-1\}
\rangle$.
\end{enumerate}
Here
$z_{\emptyset}=y_{\emptyset}:=1$,
$z_{J_1}:=z_{J_1^{(1)}}\cdots z_{J_1^{(a)}}$ and
$y_{K_1}:=y_{K_1^{(1)}}\cdots y_{K_1^{(b)}}$ 
when 
$J_1=J_1^{(1)}\cdots J_1^{(a)}$,
$K_1=K_1^{(1)}\cdots K_1^{(b)}$,
$|J_1^{(j)}|=m$ and $|K_1^{(j')}|=l$ 
for $j=1,\ldots,a$ and $j'=1,\ldots,b$.
%(ii)
\item
For 
${\cal W}_{z,y}$ in (i)(c),
${\cal W}_{z,y}\subset {\cal X}:={\rm Lin}\langle
\{\omega_{z,y}(s_{M}):
0\leq |M|\leq m-1\}
\rangle$.
\end{enumerate}
\end{lem}
%
% Proof
%
\pr
(i)
From (\ref{eqn:definitioneq}) and (\ref{eqn:omegazyb}), 
$\omega_{z,y}(X)=\omega_{z,y}(s(z)^*\,X\,s(y))$
for any $X\in \con$.
From this, 
%
% Equation 2.11
%
\begin{equation}
\label{eqn:zysj}
\omega_{z,y}(s_{J}s_{K}^*)=\overline{z_{J_1}}\,y_{K_1}
\omega_{z,y}(s_{J_2}s_{K_2}^*).
\end{equation}
%
%for any $X\in \con$, $J\in\{1,\ldots,n\}^m$
%and $K\in\{1,\ldots,n\}^l$.
If $(i,i')= (0,0)$, then (a) holds from (\ref{eqn:zysj}).
Assume $(i,i')\ne (0,0)$.
From (\ref{eqn:definitioneq}) and (\ref{eqn:omegazyb}),
%
% Equation 2.12
%
\begin{equation}
\label{eqn:abjtwo}
\begin{array}{rl}
\omega_{z,y}(s_{J_2}s_{K_2}^*)
=&
\omega_{z,y}((s(z)^*)^{\alpha}\,s_{J_2}s_{K_2}^*\,s(y)^{\beta})\\
\\
=&
\disp{\sum_{|A|=\alpha m-i}
\,
\sum_{|B|=\beta l-i'}
\overline{(z^{\otimes \alpha})_{J_2 A}}\,
(y^{\otimes \beta})_{K_2 B}
\,\omega_{z,y}(s_A^* s_B)}\\
\end{array}
\end{equation}
for any $\alpha,\beta\geq 1$.
From this, 
if $\alpha m-i=\beta l-i'$,
then 
$s_A^*s_B=\delta_{AB}I$ in (\ref{eqn:abjtwo}).
Hence (b) holds.
If $\alpha m-i\ne \beta l-i'$
for any $\alpha,\beta \geq 1$, then
%
% Equation 2.13
%
\begin{equation}
\label{eqn:sjko}
\omega_{z,y}(s_{J_2}s_{K_2}^*)
=
\left\{
\begin{array}{ll}
\disp{\sum_{|L|=l-i'}
y_{K_2 L}
\omega_{z,y}(s_{J_2}s_L)
}\quad &
\mbox{when }m-i>l-i',\\
\\
\disp{\sum_{|L|=m-i}
\overline{z_{J_2 L}}\,
\omega_{z,y}(s_L^*s_{K_2}^*)
}\quad &
\mbox{when }m-i<l-i'.\\
\end{array}
\right.
\end{equation}
If $m-i>l-i'$, then $1\leq |J_2|+|L|=i+l-i' <m$.
If $m-i<l-i'$, then $1\leq |K_2|+|L|=i'+m-i <l$.
From these and (\ref{eqn:zysj}),  (c) holds.

\noindent
(ii)
Assume $1\leq |K|\leq l-1$.
Then $l-|K|=\gamma m+j$ for some $\gamma\geq 0$ and $0\leq j\leq m-1$.
From (\ref{eqn:definitioneq}),
%
% Equation 2.14
%
\begin{equation}
\label{eqn:kmj}
\omega_{z,y}(s_K^*)
=\omega_{z,y}((s(z)^*)^{\gamma}\,s_K^*\,s(y))
=
\sum_{|J_1|=\gamma m}
\sum_{|J_2|=j}
\overline{(z^{\otimes \gamma})_{J_1}}\,y_{KJ_1J_2}\, 
\omega_{z,y}(s_{J_2}).
\end{equation}
Hence $\omega_{z,y}(s_K^*)\in {\cal X}$ and
${\cal W}_{z,y}\subset {\cal X}$.
\qedh

%
% Lemma 2.8
%
\begin{lem} 
\label{lem:correlationeq}
Assume Assumption E
for $z\in ({\cal V}_{n,m})_1$ and  $y\in ({\cal V}_{n,l})_1$
with $\omega_{z,y}$ in  (\ref{eqn:omegazyb}).
\begin{enumerate}
%(i)
\item
Assume that $d,\alpha,\beta\geq 1$ satisfy $d=\alpha m=\beta l$
and $d\geq 2$.
If $m>l$ and $z$ is nonperiodic,
then $\omega_{z,y}(s_{J})= 0$ for any $1\leq |J|\leq d-1$.
%(ii)
\item
If $m=l\geq 2$ and $z\not \sim y$, 
then $\omega_{z,y}(s_{J})= 0$ for any $1\leq |J|\leq m-1$.
%(iii)
\item
Let $\omega$ be a sub-Cuntz state by $z$ and $m\geq 2$.
If $z$ is nonperiodic,
then $\omega(s_{J})= 0$ for any $1\leq |J|\leq m-1$.
\end{enumerate}
\end{lem}
%
% Proof
%
\pr
Assume that $d,\alpha,\beta\geq 1$ satisfy $d=\alpha m=\beta l$
and $d\geq 2$.
For $1\leq a\leq d-1$,
let $J\in \{1,\ldots,n\}^a$.
From (\ref{eqn:definitioneq}),
%
% Equation 2.15
%
\begin{equation}
\label{eqn:zysjeb}
\omega_{z,y}(s_J)
=
\omega_{z,y}((s(z)^{\alpha})^*\, s_J\,s(y)^{\beta})
=
\sum_{|K|=d-a}
\sum_{|L|=a}\overline{(z^{\otimes \alpha})_{JK}}\, 
(y^{\otimes \beta})_{KL}\omega_{z,y}(s_L).
\end{equation}
By rewriting this,
%
% Equation 2.16
%
\begin{equation}
\label{eqn:doubleb}
u_a=T_{d-a}(\overline{z^{\otimes\alpha}})T_a(y^{\otimes \beta})u_a
\end{equation}
where
$T_a(z)$ is as in (\ref{eqn:ta}) and
$u_a:=\sum_{|J|=a}\omega_{z,y}(s_J)e_J\in {\cal V}_{n,a}$.

If $\omega_{z,y}(s_{J})\ne 0$ for some $J\in \{1,\ldots,n\}^a$
and $1\leq a\leq d-1$,
then $u_a\ne 0$
in (\ref{eqn:doubleb}).
Applying Lemma \ref{lem:newone}(i) to (\ref{eqn:doubleb})
with $(m,X,Y,c,v)=(d,z^{\otimes \alpha},y^{\otimes \beta},1,u_a)$,
we obtain
%
% Equation 2.17
%
\begin{equation}
\label{eqn:zoa}
z^{\otimes \alpha}=z_1\otimes z_2,\quad y^{\otimes \beta}=z_2\otimes z_1
\end{equation}
for  some $z_1\in ({\cal V}_{n,a})_1$
and $z_2\in ({\cal V}_{n,d-a})_1$.

\noindent
(i)
If $\omega_{z,y}(s_{J})\ne 0$ for some $J\in \{1,\ldots,n\}^a$
and $1\leq a\leq d-1$,
we obtain (\ref{eqn:zoa}).
From Corollary \ref{cor:projective}(iii) and $m>l$, $z$ is periodic.
This contradicts with the assumption of $z$.
Hence $\omega_{z,y}(s_{J})= 0$ for any $1\leq |J|\leq d-1$.

\noindent
(ii)
Assume $m=l$.
Then $\alpha=\beta=1$ and $d=m$ in (\ref{eqn:doubleb}).
If $\omega_{z,y}(s_{J})\ne 0$ for some $J\in \{1,\ldots,n\}^{a}$
and $1\leq a\leq m-1$,
then $z=z_1\otimes z_2$ and $y=z_2\otimes z_1$
for  some $z_1\in ({\cal V}_{n,a})_1$
and $z_2\in ({\cal V}_{n,m-a})_1$
from  (\ref{eqn:zoa}) with $\alpha=\beta=1$.
Hence $z\sim y$.
This contradicts with the assumption of $z$ and $y$.
Hence $\omega_{z,y}(s_{J})= 0$ for any $1\leq |J|\leq m-1$.

\noindent
(iii)
Remark that Assumption EC holds for $z$ and $z$.
In this case, $\omega=\omega_{z,z}$.
Hence 
(\ref{eqn:doubleb}) holds for $z=y$, $\alpha=\beta=1$ and $d=m$.
If $\omega(s_{J})\ne 0$, then 
$z_1\otimes z_2=z=z_2\otimes z_1$
for  some $z_1\in ({\cal V}_{n,a})_1$
and $z_2\in ({\cal V}_{n,m-a})_1$
from  (\ref{eqn:zoa}) with $y=z$.
From Corollary \ref{cor:projective}(i) with $c=1$,
$z$ is periodic.
This contradicts with the assumption of $z$.
Hence $\omega(s_{J})= 0$ for any $1\leq |J|\leq m-1$.
\qedh

%%%%%%%%%%%%%%%%%%%%%%%%%%%%
%
% Section 3
%
\sftt{Proofs of main theorems}
\label{section:third}
%
%%%%%%%%%%%%%%%%%%%%%%%%%%%%%%
%
% subsection 3.1
%
\ssft{Proof of Theorem \ref{Thm:maintwo}}
\label{subsection:thirdone}
%{\it Proof of Theorem \ref{Thm:maintwo}.}
When $m=1$, $\omega$ is a Cuntz state. 
Hence it suffices to show the case of $m\geq 2$.
For $J,K\in\bigcup_{a\geq 0}\{1,\ldots,n\}^a$,
we compute the value $\omega(s_Js_K^*)$ as follows.

Assume $|J|-|K|\in m{\Bbb Z}$.
Then either $|J|,|K|\in m{\Bbb Z}_{\geq 0}$ or 
$|J|,|K|\not \in m{\Bbb Z}_{\geq 0}$ holds.
If $|J|,|K|\in m{\Bbb Z}_{\geq 0}$, 
then  $\omega(s_Js_K^*)=\overline{z_J}\,z_K$
from Theorem \ref{Thm:bj}(ii)
where we use the notation in 
Theorem \ref{Thm:maintwo}(iii).
If $|J|,|K|\not \in m{\Bbb Z}_{\geq 0}$,
then  
the condition for $J,K$ in  ``otherwise"
of Theorem \ref{Thm:maintwo}(iii) holds.
In this case,
%
% Equation 3.1
%
\begin{equation}
\label{eqn:induction}
\begin{array}{rl}
\omega(s_{J}s_K^*)
=&\overline{z_{J_1}}z_{K_1}\omega(s_{J_2}s_{K_2}^*)\\
\\
=&\disp{\overline{z_{J_1}}\,z_{K_1}\sum_{|L|=m-|J_{2}|}
\omega(s_{J_2}s_Ls_L^*s_{K_2}^*)}\\
\\
=&\disp{\overline{z_{J_1}}\,z_{K_1}\sum_{|L|=m-|J_{2}|}
\overline{z_{J_2L}}\,z_{K_2L}.}\\
\end{array}
\end{equation}

Assume $|J|-|K|\not\in m{\Bbb Z}$. 
In Lemma \ref{lem:redthree},
let $l=m$ and $y=z$. Then 
the GNS representation $({\cal H},\pi,\Omega)$ by $\omega$
satisfies Assumption EC for $\Omega_z=\Omega_y=\Omega$ with
$\omega=\omega_{z,y}$. 
Applying Lemma \ref{lem:redthree}(i)(c) to $\omega$,
%
% Equation 3.2
%
\begin{equation}
\label{eqn:zzomega}
\omega(s_Js_{K}^*)\in {\cal W}_{z,z}={\rm Lin}\langle
\{\omega(s_{J}),\omega(s_K^*):1\leq |J|,|K|\leq m-1\}
\rangle.
\end{equation}

\noindent
(i)
Assume that $z$ is nonperiodic.
Then $\omega(s_J)= 0$ for all $1\leq |J|\leq m-1$
from Lemma \ref{lem:correlationeq}(iii).
From (\ref{eqn:zzomega}),
$\omega(s_Js_{K}^*)=0$ for all $J,K$ which satisfy $|J|-|K|\not\in m{\Bbb Z}$. 
From this and the case of $|J|-|K|\in m{\Bbb Z}$,
$\omega(s_{J}s_{K}^*)$ is determined by only $z$
for all $J,K$.
Hence $\omega$ is unique.

Assume $z=x^{\otimes p}$ for some $p\geq 2$ and $x\in ({\cal V}_{n,m'})_1$.
Let $\zeta:=e^{2\pi\sqrt{-1}/p}\in U(1)$.
For $1\leq j\leq p$,
assume that $\omega_j$ is a sub-Cuntz state on $\con$ by $\zeta^j x$.
Then we see that $\omega_j$ is a sub-Cuntz state by $z$ for all $j$.
On the other hand,
if $i\ne j$, then $\omega_i\ne \omega_j$.
Hence a sub-Cuntz state by $z$
is not unique.

\noindent
(ii)
From (i) and Proposition \ref{prop:stateext}(ii),
the statement holds.

\noindent
(iii)
From the proof of (i), the statement holds.
\qedh

%%%%%%%%%%%%%%%%%%%%%%%
%
% subsection 3.2
%
\ssft{Proof of  Theorem \ref{Thm:periodic}}
\label{subsection:thirdtwo}
%
%In this subsection,
%we prove Theorem \ref{Thm:periodic}.
%
% Lemma 3.2
% 
\begin{lem}
\label{lem:xvnm}
Let $x=\sum  x_Je_J\in ({\cal V}_{n,m})_1$ and $p\geq 2$.
Assume that $\con$ acts on a Hilbert space ${\cal H}$ with
a cyclic unit vector $\Omega$ which satisfies
%
% Equation 3.2
%
\begin{equation}
\label{eqn:perp}
s(x)^p\,\Omega=\Omega
\end{equation}
where $s(x):=\sum x_{J}s_J\in \con$.
Then the following hold:
\begin{enumerate}
%(i)
\item
There exist a unit vector  $(\alpha_i)\in  {\Bbb R}^p$ and 
 an orthogonal set $\{\Omega_i\}_{i=1}^{p}\subset {\cal H}$
such that 
%
% Equation 3.3
%
\begin{equation}
\label{eqn:bbb}
\Omega=\sum_{i=1}^{p}\alpha_i \Omega_i,\quad 
s(x)\Omega_i=\zeta^{-i} \Omega_i,\quad
\|\Omega_i\|=0 \mbox{ or }1
\quad (i=1,\ldots,p)
\end{equation}
where $\zeta:=e^{2\pi\sqrt{-1}/p}$.
%(ii)
\item
In (i),
define $\omega_i:=\langle\Omega_i|(\cdot)\Omega_i\rangle$.
If $\Omega_i\ne 0$,
then $\omega_i(s(x)^{k})=\zeta^{-ik}$
for any $k\geq 1$.
%(iii)
\item
Fix $i,j\in\{1,\ldots,p\}$ and 
define $x':=\zeta^{i}x$ and  $x'':=\zeta^{j}x$.
Assume $m\geq 2$.
If $i\ne j$, then 
for $J,K\in\bigcup_{a=1}^{m-1}\{1,\ldots,n\}^a$,
%
% Equation 3.4
%
\begin{equation}
\label{eqn:langleomega}
\langle\Omega_i|s_{J}s_K^*\Omega_j\rangle 
=
\left\{
\begin{array}{ll}
\disp{\sum_{|L|=m-|J|}\overline{x'_{JL}}\,
\langle \Omega_i|s_L^*s_K^*\Omega_j\rangle}\quad &
\mbox{when }|J|>|K|,\\
\\
0 \quad & \mbox{when }|J|=|K|,\\
\\
\disp{\sum_{|L|=m-|K|}x''_{KL}
\langle \Omega_i|s_Js_L\Omega_j\rangle}\quad &
\mbox{when }|J|<|K|.\\
\end{array}
\right.
\end{equation}
%(iv)
\item
If $i\ne j$, then 
$\langle\Omega_i|s_{J}s_K^*\Omega_j\rangle=0$
for any $J,K$.
\end{enumerate}
\end{lem}
%
% Proof
%
\pr
(i)
Define ${\cal K}:=\{w\in {\cal H}: s(x)^p\,w=w\}$.
Since $\Omega\in {\cal K}$,
${\cal K}$ is a non-zero closed subspace of ${\cal H}$.
Then $R:=s(x)|_{{\cal K}}$ satisfies $R^p=I_{{\cal K}}$.
From this,
we obtain the spectral decomposition 
$R=\sum_{j=1}^{p}\zeta^{-j}E_j$
where $\{E_j\}$ is the orthogonal set of projections on ${\cal K}$
such that $E_1+\cdots+E_p=I_{{\cal K}}$.
For $i=1,\ldots,p$,
define $\alpha_i:=\|E_i\Omega\|$, 
and $\Omega_i:=\alpha_i^{-1}E_i\Omega$ when $\alpha_i\ne 0$ and
$\Omega_i:=0$ otherwise.
Then the statement holds.

\noindent
(ii) From (\ref{eqn:bbb}),
the statement holds.

\noindent
(iii)
From (\ref{eqn:bbb}), 
%
% Equation 3.5
%
\begin{equation}
\label{eqn:xxdash}
s(x')\Omega_i=\Omega_i,\quad s(x'')\Omega_j=\Omega_j.
\end{equation}
In Lemma \ref{lem:redthree},
let $l=m$, $z:=x'$ and $y:=x''$.
Assumption E holds for $z$ and $y$ with 
$\omega_{z,y}=\langle \Omega_i|(\cdot)\Omega_j\rangle$.
Assume $|J|=|K|=a$.
From (i),
$\omega_{z,y}(I)=\langle\Omega_i|\Omega_j\rangle=0$.
From this, the statement holds.
The rest  is proved from Lemma \ref{lem:redthree}(i).

\noindent
(iv)
If $m=1$,
then $\langle\Omega_i|s_Js_K^*\Omega_j\rangle
=\overline{x'_J}\,x''_K\langle\Omega_i|\Omega_j\rangle=0$
for any $J,K$ when $i\ne j$.
Hence we assume $m\geq 2$.
From (iii) and Lemma \ref{lem:redthree}(ii),
it is sufficient to show 
$\langle\Omega_i|s_{J}\Omega_j\rangle=0$
for all $1\leq |J|\leq m-1$.
In Lemma \ref{lem:correlationeq}(ii),
let $z:=x'$ and $y:=x''$.
Since $z\not \sim y$,
$\langle\Omega_i|s_{J}\Omega_j\rangle=0$
for $1\leq |J|\leq m-1$
from Lemma \ref{lem:correlationeq}(ii).
\qedh

\ww
{\it Proof of Theorem \ref{Thm:periodic}.}
Assume that $\omega$ is a  sub-Cuntz state by $z$
with the GNS representation $({\cal H},\pi,\Omega)$.
For $X\in \con$,
we write $X$ as $\pi(X)$ for the simplicity of description.
Then the assumption in Lemma \ref{lem:xvnm} is satisfied.
For $\alpha_i$ in Lemma \ref{lem:xvnm}(i),
define $a_i:=\alpha_i^2$ for $i=1,\ldots,p$.
From Lemma \ref{lem:xvnm}(i) and (iv),
$\omega
=\sum_{i=1}^p|\alpha_i|^2
\omega_{i}
=\sum_{i=1}^pa_i\omega_{i}$.
If $\Omega_j\ne 0$, then 
we see that
$\omega_j$ is the sub-Cuntz state by $\zeta^j x$.
If $\Omega_j=0$, then
$\omega_j=0$ and $a_j=0$.
In this case, we can replace $\omega_j$ with the sub-Cuntz state
by $\zeta^j x$ with keeping $\sum_{i=1}^pa_i\omega_{i}$.
Therefore (\ref{eqn:omegatwo}) holds
as a convex-hull of states.

We prove the uniqueness as follows.
Assume $\omega=\sum_{j=1}^{p}b_j\omega_j$
for some $(b_1,\ldots,b_p)\in\Delta_{p-1}$.
From Lemma \ref{lem:xvnm}(ii),
%
% Equation 3.10
%
\begin{equation}
\label{eqn:minus}
0=\sum_{j=1}^{p}(a_j-b_j)\omega_j(s(x)^k)
=\sum_{j=1}^{p}(a_j-b_j)\zeta^{-jk}
\quad \mbox{ for all } k\geq 1.
\end{equation}
This implies $a_j-b_j=0$ for all $j$.
Hence $(a_1,\ldots,a_p)$  is unique.
\qedh

%%%%%%%%%%%%%%%%%%%%%
%
% subsection 3.3
%
\ssft{Proof of  Theorem \ref{Thm:mainthree}}
\label{subsection:thirdthree}
%
%In this subsection, we prove Theorem \ref{Thm:mainthree}.
Let $\sim$ be as in Theorem \ref{Thm:mainthree}.
By Theorem \ref{Thm:bj},  we can prove the following.
%
% Lemma 3.3
%
\begin{lem}
\label{lem:asse}
For two nonperiodic parameters
$z\in ({\cal V}_{n,m})_1$ and $y\in ({\cal V}_{n,l})_1$,
the following are equivalent:
\begin{enumerate}
%(i)
\item
$\tilde{\omega}_z\sim \tilde{\omega}_y$.
%(ii)
\item
Assumption EC in Lemma \ref{lem:redthree}
holds for $z$ and $y$.
\end{enumerate}
\end{lem}

%
% Lemma 3.4
%
\begin{lem}
\label{lem:vanish}
Let  $z\in ({\cal V}_{n,m})_1$ and  $y\in ({\cal V}_{n,l})_1$.
\begin{enumerate}
%(i)
\item
Assume Assumption EC for $z$ and  $y$ with 
$\omega_{z,y}$ in (\ref{eqn:omegazyb}).
If both $z$ and $y$ are nonperiodic,
then  $\omega_{z,y}\not \equiv 0$.
%(ii)
\item
Assume Assumption E for $z$ and $y$.
If there exist  integers $\alpha,\beta\geq 1$ such that $\alpha m=\beta l$
and $z^{\otimes \alpha }\ne y^{\otimes \beta}$,
then  $\omega_{z,y}(I)=0$.
Especially,
if $m=l$ and 
 $z\ne y$, then $\omega_{z,y}(I)=0$.
\end{enumerate}
\end{lem}
%
% Proof
%
\pr
(i)
Let ${\cal H}$ be as in Assumption E.
%By Lemma \ref{lem:asse},
%$\tilde{\omega}_z\sim \tilde{\omega}_y$.
If  $\omega_{z,y}\equiv 0$, then
$0=\omega_{z,y}(s_J)=\langle \Omega_z|s_J\Omega_y\rangle$
for all $J$.
Since $\langle \Omega_y|(\cdot)\Omega_y\rangle=\tilde{\omega}_y$
and 
${\cal H}$ is generated by $\{s_J\Omega_y:J\}$
from Lemma \ref{lem:omegasub}(ii), $\Omega_z=0$.
This contradicts with $\|\Omega_z\|=1$ in Assumption EC.
Hence  $\omega_{z,y}\not \equiv 0$.

\noindent
(ii)
By (\ref{eqn:definitioneq}) and (\ref{eqn:omegazyb}),
%
% Equation 3.7
%
\begin{equation}
\label{eqn:identity}
\omega_{z,y}(I)
=
\omega_{z,y}(\, (s(z)^*)^{\alpha }\, s(y)^{\beta}\,)
=\langle z^{\otimes \alpha }|y^{\otimes \beta}\rangle\omega_{z,y}(I).
\end{equation}
By Lemma \ref{lem:trivial},
$z^{\otimes \alpha }\ne y^{\otimes \beta}$
implies 
$\langle z^{\otimes \alpha}|y^{\otimes \beta}\rangle\ne 1$.
Hence $\omega_{z,y}(I)=0$.
\qedh

%
% Lemma 3.5
%
\begin{lem}
\label{lem:nes}
Let $z\in ({\cal V}_{n,m})_1$ and  $y\in ({\cal V}_{n,l})_1$.
If $l=m$ and both $z$ and $y$ are nonperiodic, then 
the following are equivalent:
\begin{enumerate}
%(i)
\item
 $\tilde{\omega}_z\sim \tilde{\omega}_y$.
%(ii)
\item
 $z\sim y$.
\end{enumerate}
\end{lem}
%
% Proof
%
\pr
When $m=1$, both $\tilde{\omega}_z$
and $\tilde{\omega}_y$ are Cuntz states. 
Hence it is sufficient to show the case of $m\geq 2$.

\noindent
(i)$\Rightarrow$(ii)
Assume  $\tilde{\omega}_z\sim \tilde{\omega}_y$.
If $z=y$, then the statement holds.
Assume  $z\ne y$.
From Lemma \ref{lem:asse}, Assumption EC for $z$ and $y$ holds.
Let  $\omega_{z,y}$ be as in (\ref{eqn:omegazyb}).
Then $\omega_{z,y}(I)=0$
from Lemma \ref{lem:vanish}(ii).
From Lemma \ref{lem:vanish}(i) and Lemma \ref{lem:redthree}(ii),
there must exist $1\leq a\leq m-1$ such that 
$\omega_{z,y}(s_J)\ne 0$ for some $J\in \{1,\ldots,n\}^a$.
From Lemma \ref{lem:correlationeq}(ii), $z\sim y$.

\noindent
(ii)$\Rightarrow$(i)
Assume $z\sim y$ and $z\ne y$.
Then there exist $x_1,x_2\in 
\bigcup_{a\geq 1}({\cal V}_{n,a})_1$ such that
$z=x_1\otimes x_2$ and $y=x_2\otimes x_1$.
Let $({\cal H},\pi,\Omega)$ denote
the GNS representation by $\tilde{\omega}_z$.
From Theorem \ref{Thm:bj}(ii),
 $\pi(s(z))\Omega=\Omega$.
Then $\Omega':=\pi(s(x_2))\Omega\in {\cal H}$ 
is also a cyclic unit vector because $\tilde{\omega}_z$ is pure, and
we can verify  $\pi(s(y))\Omega'=\Omega'$.
Hence Assumption EC for $z$ and $y$ holds.
From Lemma \ref{lem:asse}, 
 $\tilde{\omega}_z\sim\tilde{\omega}_y$.
\qedh

\ww
{\it Proof of Theorem \ref{Thm:mainthree}.}
(i)$\Rightarrow$(ii) 
Assume $\tilde{\omega}_z\sim\tilde{\omega}_y$.
From Lemma \ref{lem:asse},
Assumption EC for $z$ and $y$ holds.
Let $\omega_{z,y}$ be as in (\ref{eqn:omegazyb}).

Assume $m>l$.
From Corollary \ref{cor:projective}(ii),
there exist no $\alpha,\beta$ such that
$z^{\otimes \alpha}=y^{\otimes \beta}$.
From this and Lemma \ref{lem:vanish}(ii),
$\omega_{z,y}(I)=0$.
Since $z$ is nonperiodic, 
$\omega_{z,y}(s_J)=0$ for any $1\leq |J|\leq d-1$
from Lemma \ref{lem:correlationeq}(i).
From Lemma \ref{lem:redthree}(i) and (ii),
$\omega_{z,y}\equiv 0$.
From Lemma \ref{lem:vanish}(i),
Assumption EC does not hold.
From Lemma \ref{lem:asse},
$\tilde{\omega}_z\not\sim \tilde{\omega}_y$.
Hence $m\not >l$.
As the same token, 
we obtain $m\not <l$.
Hence $m=l$.
From Lemma \ref{lem:nes}, $z\sim y$.

\noindent
(ii)$\Rightarrow$(i) 
Assume $z\sim y$. Then $m=l$.
From Lemma \ref{lem:nes},
$\tilde{\omega}_z\sim\tilde{\omega}_y$. 
\qedh

%%%%%%%%%%%%%%%%%%%%
%
% subsection 3.4
%
\ssft{Proofs of Proposition \ref{prop:quasitensor}\
and Proposition \ref{prop:tensor}}
\label{subsection:thirdfour}

%\ww
\noindent
{\it Proof of Proposition \ref{prop:quasitensor}.}
Let $\omega$ and $\omega'$ 
be sub-Cuntz states by $z\in ({\cal V}_{n,m})_1$
and $y\in ({\cal V}_{n',l})_1$, respectively.
Then $z^{\otimes l}\in ({\cal V}_{n,ml})_1$
and $y^{\otimes m}\in ({\cal V}_{n',ml})_1$.
For any $J\in\{1,\ldots,nn'\}^{ml}$,
we can verify
$(\omega\ptimes\omega')(s_{J}^{(nn')})
=
\overline{(z^{\otimes l}\boxtimes  y^{\otimes m})_{J}}$
where $\boxtimes$ is as in (\ref{eqn:boxtimes}).
Hence 
$\omega\ptimes\omega'$ is a sub-Cuntz state by
$z^{\otimes l}\boxtimes  y^{\otimes m}\in ({\cal V}_{nn',\,ml})_1$.
\qedh

\noindent
In this proof,
there is no assumption of nonperiodicity for $z$ and $y$.
Hence $\omega,\omega'$ and $\omega\ptimes \omega'$
are not always unique.
%
% Lemma 3.5
%
\begin{lem}
\label{lem:boxtimes}
\begin{enumerate}
%(i)
\item
For $x,x'\in{\cal V}_{n,m}$ and $y,y'\in {\cal V}_{n',m}$,
$\langle x\boxtimes y|x'\boxtimes y'\rangle =
\langle x|x'\rangle 
\langle y|y'\rangle$.
%(ii)
\item
Let $x,x'\in ({\cal V}_{n,m})_1$ and $y,z\in ({\cal V}_{n',m})_1$.
If  $x\boxtimes y=x'\boxtimes z$
or $y\boxtimes x=z\boxtimes x'$,
then $y=cz$ for some $c\in U(1)$. 
%(iii)
\item
If  $x\in ({\cal V}_{n,m})_1$
and  $y\in ({\cal V}_{n',m})_1$
satisfy $x\boxtimes y= w^{\otimes p}$
for some $w\in ({\cal V}_{nn',m'})_1$
and $p\geq 2$,
then 
there exist  $v_1\in ({\cal V}_{n,m'})_1$
and  $v_2\in ({\cal V}_{n',m'})_1$
such that  $x=v_1^{\otimes p}$
and $y=v_2^{\otimes p}$.
\end{enumerate}
\end{lem}
%
% Proof
%
\pr
(i)
By definition, the statement holds.

\noindent
(ii)
Assume $x\boxtimes y=x'\boxtimes z$.
By assumption and (i),
$1
=\langle 
x\boxtimes y|x'\boxtimes z
\rangle 
=
\langle 
x|x'\rangle 
\langle 
y|z\rangle$.
By applying Lemma \ref{lem:trivial} to this,
the statement holds.
As the same token, the rest is proved. 

\noindent
(iii)
By assumption, $m'p=m$.
For $J\in \{1,\ldots,nn'\}^{m'}$,
$(x\boxtimes y)_{J^p}=(w_J)^p$
where $J^p=J\cdots J$ ($p$-times).
Hence
$x_{J_1^p} y_{J_2^p}=(w_{J_1\boxtimes J_2})^p$
for all
 $J_1\in \{1,\ldots,n\}^{m'}$
and 
 $J_2\in \{1,\ldots,n'\}^{m'}$.
Then there exist 
$p$-th roots  $A_{J_1}$ and  $B_{J_2}$
of $x_{J_1^p}$ and  $y_{J_2^p}$
such that
$A_{J_1}B_{J_2}=w_{J_1\boxtimes J_2}$.
Define
$A:=\sum_{J_1}A_{J_1}e_{J_1}\in {\cal V}_{n,m'}$ and 
$B:=\sum_{J_2}B_{J_2}e_{J_2}\in {\cal V}_{n',m'}$.
Then $w=A\boxtimes B$.
By normalizing $A$ and $B$,
we obtain two unit vectors $w_1,w_2$ such that 
$w=w_1\boxtimes w_2$.
% and $\|w_1\|=\|w_2\|=1$.
From these,
$x\boxtimes y=w^{\otimes p}
=(w_1\boxtimes w_2)^{\otimes p}
=w_1^{\otimes p}\boxtimes w_2^{\otimes p}$.
From (ii), 
$x=cw_1^{\otimes p}$
and 
$y=\overline{c}w_2^{\otimes p}$
for some $c\in U(1)$.
From these,
we can choose $v_1$ and $v_2$ as the statement.
\qedh

\ww
{\it Proof of Proposition \ref{prop:tensor}.}
Let $\alpha,\beta,d$ be as in (\ref{eqn:albe}).

\noindent
(i)
Assume  $z*y=w^{\otimes p}$ for some 
$w\in {\cal V}_{nn',k}$ and $p\geq 2$
where $k:=d/p$.
By definition,
$z^{\otimes \alpha}\boxtimes y^{\otimes \beta}=w^{\otimes p}$.
From Lemma \ref{lem:boxtimes}(iii),
we obtain $v_1\in ({\cal V}_{n,k})_1$ and $v_2\in ({\cal V}_{n',k})_1$ 
such that
%
% Equation 3.9
%
\begin{equation}
\label{eqn:vone}
z^{\otimes \alpha}=v_1^{\otimes p},\quad 
y^{\otimes \beta}=v_2^{\otimes p}.
\end{equation}
From Corollary \ref{cor:projective}(ii),
$v_1=c'' z^{\otimes d_1}$
for some $d_1\geq 1$ and $c''\in U(1)$.
From this and (\ref{eqn:vone}),
$z^{\otimes \alpha}=
v_1^{\otimes p}=(c''z^{\otimes d_1})^{\otimes p}
=(c'')^pz^{\otimes d_1p}$.
Hence $\alpha =p d_1$.
As the same token,
$\beta=pd_2$ for some $d_2\geq 1$.
Therefore 
 $\alpha$ and $\beta$ have a common divisor $p\geq 2$.
This contradicts with the choice of $\alpha$ and $\beta$.
Therefore $z*y$ is nonperiodic.

\noindent
(ii)
Remark that $\tilde{\omega}_{z* y}$ is uniquely defined by (i)
and Theorem \ref{Thm:maintwo}(i).
We see that
$\{\tilde{\omega}_z\ptimes \tilde{\omega}_y\}(s_J^{(nn')})
=\overline{(z* y)_{J}}$
for all $J\in\{1,\ldots,nn'\}^d$.
Hence 
$\tilde{\omega}_z\ptimes \tilde{\omega}_y$ is a sub-Cuntz state
by $z* y$.
Since 
a sub-Cuntz state by $z* y$ is unique,
the statement holds.
\qedh

%%%%%%%%%%%%%%%%%%%%%%%%%%%%%
%
% Section 4
%
\sftt{Examples}
\label{section:fourth}
In this section, 
we show examples so that a reader can easily check main theorems.
%%%%%%%%%%%%%%
%
% subsection 4.1
%
\ssft{Sub-Cuntz states of order $2$}
\label{subsection:fourthone}
In this subsection, we show sub-Cuntz states on $\con$ of order $2$ 
as simplest, nontrivial and essentially new examples of main theorems.
For convenience,
we rewrite main theorems 
in $\S$ \ref{subsection:firsttwo}
for the case of $m=2$ as follows.
%
% Theorem 4.1
% 
\begin{Thm}
\label{Thm:mainone}
Let $(({\Bbb C}^n)^{\otimes 2})_1:
=\{z\in {\Bbb C}^n\otimes {\Bbb C}^n:\|z\|=1\}$.
Fix $z=\sum_{ij}z_{ij}e_{i}\otimes e_{j}
\in (({\Bbb C}^n)^{\otimes 2})_1$.
Let $\omega$ be a sub-Cuntz state on $\con$  by $z$, that is,
$\omega$ is a state on $\con$ which satisfies
%
% Equation 4.1
%
\begin{equation}
\label{eqn:mtwo}
\omega(s_{i}s_{j})=\overline{z_{ij}}\quad \mbox{for all }i,j=1,\ldots,n.
\end{equation}
Then such $\omega$ exists for any $z$ and the following hold:
\begin{enumerate}
%(i)
\item
$\omega$ is unique if and only if 
$z$ is nonperiodic, that is $z\not\in \{x\otimes x:x\in {\Bbb C}^n\}$.
In this case,
we write $\tilde{\omega}_z$ as $\omega$.
%(ii)
\item
If $z$ is nonperiodic, then
$\tilde{\omega}_z$ is pure, and the following holds:
%
% Equation 4.2
%
\begin{equation}
\label{eqn:omegasjskb}
\tilde{\omega}_z(s_Js_{K}^*)=\left\{
\begin{array}{ll}
\overline{z_{J}}\,z_{K}\quad &\mbox{when both $|J|$ and $|K|$ are even},\\
\\
\overline{z_{J_1}}\,z_{K_1}
\disp{\sum_{d=1}^n\overline{z_{jd}}\,z_{kd}}
\quad &\mbox{when }
\begin{array}{l}
 J=J_1j,\,K=K_1k,\\
\mbox{both $|J_1|$ and  $|K_1|$ are even},\\
\end{array}
\\
\\
0 \quad& \mbox{when } |J|-|K|\mbox{ is odd}\\
\end{array}
\right.
\end{equation}
for $J,K\in\bigcup_{a\geq 1}\{1,\ldots,n\}^a\cup\{\emptyset\}$
where 
$z_J:=z_{J^{(1)}}\cdots z_{J^{(l)}}$ when $J=J^{(1)}\cdots J^{(l)}$
and $|J^{(i)}|=2$ for $i=1,\ldots,l$.
%(iii)
\item
%(degenerate case)
If $z=x\otimes x$ for some $x=(x_1,\ldots,x_n)\in{\Bbb C}^n$,
then there exists a real number $0\leq a\leq 1$ such that 
$\omega$ has the following form
%
% Equation 4.3
%
\begin{equation}
\label{eqn:omeaa}
\omega=a\omega_++(1-a)\omega_-
\end{equation}
where $\omega_\pm$ denotes
the Cuntz state on $\con$ by $\pm x$, 
that is,
$\omega_{\pm}$ satisfies
$\omega_{\pm}(s_i)=\pm\overline{x_i}$ for all $i$.
%(iv)
\item
Let $z,y\in (({\Bbb C}^{n})^{\otimes 2})_1$.
If both $z$ and $y$ are nonperiodic,
then $\tilde{\omega}_z\sim\tilde{\omega}_y$
if and only if {\rm (a)} $z=y$, or {\rm (b)} $z=x_1\otimes x_2$
and $y=x_2\otimes x_1$ for some 
$x_1,x_2\in {\Bbb C}^n$.
\end{enumerate}
\end{Thm}
%
% Proof
%
\pr
The existence of $\omega$ holds from  Fact \ref{fact:existence}.

\noindent
(i)
From Theorem \ref{Thm:maintwo}(i),
the statement holds.

\noindent
(ii)
From Theorem \ref{Thm:maintwo}(ii) and (iii),
statements hold.

\noindent
(iii)
From the case of $(p,m')=(2,1)$ in Theorem \ref{Thm:periodic},
the statement holds.

\noindent
(iv)
From Theorem \ref{Thm:mainthree}, the statement holds.
\qedh

We show a more convenient corollary as follows.
%
% Corollary 4.2
%
\begin{cor}
\label{cor:nonsym}
Assume the same assumption in Theorem \ref{Thm:mainone}
for $z=\sum_{ij}z_{ij}e_i\otimes e_j$.
\begin{enumerate}
%(i)
\item
If $A:=(z_{ij})\in M_n({\Bbb C})$ satisfies
$\|A\|<1$, then $\omega$ is unique and pure.
%(ii)
\item
If $z_{ij}\ne z_{ji}$ for some $i,j$,
then $\omega$ is unique and pure.
\end{enumerate}
\end{cor}
%
% Proof
% 
\pr
(i)
Remark that $A$ coincides with $T_1(z)$ in (\ref{eqn:ta})
as operators on ${\Bbb C}^n$.
The assumption implies that $z$ is indecomposable.
Especially, $z$ is nonperiodic.
From Theorem \ref{Thm:mainone}(i) and (ii), the statement holds.

\noindent
(ii)
In this case, $z$ is nonperiodic.
Hence the statement holds
from Theorem \ref{Thm:mainone}(i) and (ii).
\qedh

Next, we show concrete examples. 
In stead of $z=\sum_{ij}z_{ij}e_{i}\otimes e_{j}
\in
(({\Bbb C}^n)^{\otimes 2})_1$,
we use a matrix $A=(z_{ij})\in M_n({\Bbb C})$ such that $\|A\|_2=1$
in order to apply Corollary \ref{cor:nonsym}.
We assume that $\con$ acts on a Hilbert space with a cyclic unit vector
$\Omega$.
Define the vector state $\omega$ on $\con$ with respect to $\Omega$:
%
% Equation 4.4
%
\begin{equation}
\label{eqn:omegaeq}
\omega:=\langle\Omega|(\cdot)\Omega\rangle.
\end{equation}
%
% Example 4.3
%
\begin{ex}
\label{ex:narray}
{\rm
Let  $(c_i)\in {\Bbb C}^n$ be a unit vector.
Assume that the following equation holds:
%
% Equation 6.5
%
\begin{equation}
\label{eqn:cisi}
\sum_{i=1}^{n}c_{i}s_{i}^2\Omega=\Omega.
\end{equation}
\begin{enumerate}
%(i)
\item
If  $|c_i|<1$ for all $i$, 
then $A:={\rm diag}(c_1,\ldots,c_n)\in M_{n}({\Bbb C})$
satisfies $\|A\|<1$.
Hence  $\omega$ is unique and pure
from 
Corollary \ref{cor:nonsym}(i).
%(ii)
\item
If there exists $i$ such that $|c_i|=1$,  then
$c_is_i^2\Omega=\Omega$.
Let $q\in U(1)$ be a quadratic root of $c_i$.
From Theorem \ref{Thm:mainone}(iii),
there exists $0\leq a\leq 1$ such that
$\omega=a\omega_++(1-a)\omega_-$
where $\omega_{\pm}$ denotes
the Cuntz state by $\pm qe_i$, that is,
$\omega_{\pm}$ satisfies
$\omega_{\pm}(s_i)=\pm\overline{q} $.
In this case,
$\omega$ is pure if and only if  $a= 0$ or $1$.
%Especially,  if $c_i=1$,
%then $\omega$ is a permutative state
%(see $\S$ \ref{subsection:fourthtwo}).
%In this case,  $\omega_{\pm}(s_i)=\pm 1$.
%Hence $\omega_{-}$ is not permutative.
\end{enumerate}
}
\end{ex}

Fix $n=2$ from here.
Let $s_{ij}:=s_is_j$ for $i,j=1,2$.
%
% Example 4.4
%
\begin{ex}
\label{ex:assume}
{\rm
Assume that the following equation holds:
%
% Equation 4.6
%
\begin{equation}
\label{eqn:frac}
\frac{1}{2}(s_{11}-s_{12}+s_{21}+s_{22})\Omega=\Omega.
\end{equation}
Then $(z_{ij})=\frac{1}{2}\left[
\begin{array}{cc}
1 &-1\\
1 &1
\end{array}
\right]$
satisfies $z_{12}\ne z_{21}$.
From Corollary \ref{cor:nonsym}(ii), $\omega$ is unique and pure.
}
\end{ex}
%
% Example 4.5
%
\begin{ex}
\label{ex:assumethat}
{\rm
Assume that the following equation holds:
%
% Equation 4.7
%
\begin{equation}
\label{eqn:fracone}
\frac{1}{\sqrt{2}}(s_{12}+s_{21})\Omega=\Omega.
\end{equation}
Then $A=\frac{1}{\sqrt{2}}\left[
\begin{array}{cc}
0 &1\\
1 &0
\end{array}
\right]$
satisfies $\|A\|=\frac{1}{\sqrt{2}}<1$.
From  Corollary \ref{cor:nonsym}(i),
$\omega$ is unique and pure.
}
\end{ex}

%%%%%%%%%%%%%%%%%%%%%%%%%%%%%
%
% subsection 4.2
%
\ssft{Sub-Cuntz states 
associated with permutative representations}
\label{subsection:fourthtwo}
In this subsection,
we show known results in $\S$ 5 of \cite{BJ,TS01}
by using results of sub-Cuntz states.
A representation $({\cal H},\pi)$ of $\con$ is said to be {\it permutative} 
if there exists an orthonormal basis $B=\{v_k:k\in\Lambda\}$ of ${\cal H}$
such that $\pi(s_i) v_k\in B$  for any $i,k$ \cite{BJ,DaPi2,DaPi3}.
We explain sub-Cuntz states associated with 
permutative representations as follows.
For $z=\sum z_Je_J\in ({\cal V}_{n,m})_1$,
assume $z_J=1$ for some $J$.
In this case, $z=e_J$ and the following holds.
%
% Proposition 4.6
%
\begin{prop}
\label{prop:permutwo}
\begin{enumerate}
%(i)
\item
For any  $J\in \{1,\ldots,n\}^m$,
there exists a state $\omega$ on $\con$ which satisfies
%
% Equation 4.8
%
\begin{equation}
\label{eqn:omeper}
\omega(s_J)=1.
\end{equation}
%
%(ii)
\item
If a state $\omega$ on $\con$ satisfies (\ref{eqn:omeper}),
then $\omega$ is unique if and only if $J$ is nonperiodic,
that is,
$J=K^p$ for some $K$ implies $p=1$.
In this case, $\omega$ is a pure sub-Cuntz state $\tilde{\omega}_{e_J}$ 
by $e_J$, and we write $\phi_J$ as $\tilde{\omega}_{e_J}$ for short.
%(iii)
\item
If $J$ is nonperiodic, then
the GNS representation by $\phi_J$ is permutative.
%(iv)
\item
If $J=K^p$ for some nonperiodic word $K$ and $p\geq 2$,
then 
$\phi_J=\sum_{j=1}^{p}a_j\omega_j$
for some nonnegative numbers $a_1,\ldots,a_p$ such that 
$a_1+\cdots+a_p=1$
where $\omega_j$ denotes the sub-Cuntz state by 
$e^{2\pi j\sqrt{-1}/p}e_J$.
%(v)
\item
For two nonperiodic words $J$ and $K$,
$\phi_J\sim \phi_K$ if and only if 
$J$ and $K$ are conjugate,
that is,
$J=K$ or $J=L_1L_2$ and $K=L_2L_1$ for some $L_1,L_2$.	
%(vi)
\item
Let ${\goth S}_n$ denote the symmetric group 
on the set $\{1,\ldots,n\}$.
Define the action of ${\goth S}_n$ on ${\Bbb C}^n$
as  $\sigma e_i:=e_{\sigma(i)}$ 
for $i=1,\ldots,n$ $\sigma\in {\goth S}_n$.
With respect to this action,
we identify ${\goth S}_n$ with the subgroup of $U(n)$.
Then
for any nonperiodic word $J$,
$\alpha^*_{\sigma}\circ \phi_J=\phi_{\sigma J}$
for any $\sigma\in {\goth S}_n$
where $\alpha^*$ is as in $\S$ \ref{subsubsection:firstthreeone}
and 
$\sigma J:=(\sigma(j_1),\ldots,\sigma(j_l))$
when $J=(j_1,\ldots,j_l)$.
%(vii)
\item
Let $\ptimes$ and $\boxtimes$ be as in $\S$ \ref{subsubsection:firstthreetwo}.
For two nonperiodic words $J$ and $K$,
$\phi_J\ptimes \phi_K=\phi_{J*K}$
where
$J*K:=J^{\alpha}\boxtimes K^{\beta}$
such that 
$\alpha,\beta\geq 1$ and
$\alpha |J|=\beta|K|$ is the least common multiple
of $|J|$ and $|K|$.
\end{enumerate}
\end{prop}
%
% Proof
%
\pr
(i) Let $z:=e_J\in ({\cal V}_{n,m})_1$ when $|J|=m$.
Then we see that $\omega$ satisfies
$\omega(s_{K})=\delta_{JK}$ for
all $K\in\{1,\ldots,n\}^m$.
Hence $\omega$ is a sub-Cuntz state by $z$.
From  Fact \ref{fact:existence}, the statement holds.

\noindent
(ii)
Let $z:=e_J\in ({\cal V}_{n,m})_1$  when $|J|=m$.
By assumption,  $z$ is nonperiodic if and only if $J$ is nonperiodic.
From Theorem \ref{Thm:maintwo}(i) and (ii),
the statement holds.

\noindent
(iii)
Let $({\cal H},\pi,\Omega)$ denote
the GNS representation by $\phi_J$.
We prove that  $({\cal H},\pi)$ is permutative as follows.
Since $\pi(s_J)\Omega=\Omega$,  
$\phi_J(s_{K})=1$ when $K\in X:=\{J^a:a\geq 0\}$
and $\phi_J(s_{K})=0$ when $K\not \in X$.
Let $v_K:=\pi(s_K)\Omega$.
Since $J$ is nonperiodic,
$\langle v_K|v_L\rangle =1$ when 
$K=LJ^a$ or $L=KJ^a$ for some $a\geq 0$
and 
$\langle v_K|v_L\rangle =0$ otherwise.
From Lemma \ref{lem:trivial},
$\langle v_K|v_L\rangle =1$ if and only if $v_K=v_L$.
Therefore
$\{u\in {\cal H}: \mbox{there exists } K \mbox{ such that }u=v_K\}$
 is an orthonormal basis of ${\cal H}$ from
Lemma \ref{lem:omegasub}(ii).
Hence $({\cal H},\pi)$ is permutative.

\noindent
(iv)
From Theorem \ref{Thm:periodic},
the statement holds.

\noindent
(v)
From Theorem \ref{Thm:mainthree},
the statement holds.

\noindent
(vi)
From Proposition \ref{prop:cova},
the statement holds.

\noindent
(vii)
From Proposition \ref{prop:tensor},
the statement holds.
\qedh

\noindent
When $m=1$ in Proposition \ref{prop:permutwo},
there exists $i\in\{1,\ldots,n\}$ such that $z=e_i$.
In this case, $z$ is nonperiodic and the following holds:
%
% Equation 4.9
%
\begin{equation}
\label{eqn:perm}
\omega(s_j)=\delta_{ij}\quad \mbox{for all }j=1,\ldots,n.
\end{equation}
Any Cuntz state is given as 
a transformation of this by 
the dual action of the standard $U(n)$-action on $\con$
(see the proof of Theorem \ref{Thm:cuntzstate}(ii)).
%
% Fact 4.7
%
\begin{fact}
%(\cite{BJ})
\label{fact:vector}
Assume that $\con$ acts on the Hilbert space 
$\ell^2(\Lambda)$ 
with an orthonormal basis $B=\{v_{\lambda}:\lambda\in\Lambda\}$
such that 
$s_iv_{\lambda}\in B$ for any $i$ and $\lambda$,
and 
$\omega$ is the vector state on $\con$ by $v_{\lambda_0}$
for some $\lambda_0\in\Lambda$.
If $\omega$ is finitely correlated, 
then $\omega_L:=\omega(s_{L}(\cdot)s_{L}^*)$ 
is a sub-Cuntz state for some $L\in\{1,\ldots,n\}^k$.
Especially,
if $\omega$ is pure, then $\omega_L$ is also pure
and $\omega_L\sim \omega$.
\end{fact}
%
% Proof
%
\pr
Let $\Omega:=v_{\lambda_0}$.
By assumption,
the action of $\con$ on ${\cal H}$ is permutative.
Hence $s_J^*\Omega=0$ or 
$s_J^*\Omega\in B$ for any $J$.
Therefore $\{s_J^*\Omega:J\}\subset B\cup \{0\}$ and
$\#\{s_J^*\Omega:J\}= \#(\{s_J^*\Omega:J\}\cap B)+1$.
By assumption,
%
% Equation 4.10
%
\begin{equation}
\label{eqn:infinitedim}
\infty>\dim{\rm Lin}\langle \{s_J^*\Omega:J\}\rangle
=\dim{\rm Lin}\langle \{s_J^*\Omega\in B:J\}\rangle
=\#\{s_J^*\Omega\in B:J\}.
\end{equation}
Hence $\#\{s_J^*\Omega:J\}<\infty$.
For any $\lambda\in\Lambda$,
there exists a unique $i$
such that $s_i^*v_{\lambda}\in B$.
Hence  there exists a unique sequence
$\{J^{(l)}\in\bigcup_{a\geq 1}\{1,\ldots,n\}^{a}:|J^{(l)}|=l\mbox{ for all }l\}$ 
such that 
$s_{J^{(l)}}^*\Omega \ne 0$
for any $l$.
Since $\#\{s_J^*\Omega:J\}<\infty$,
there exist $p\geq 1$ 
and $l_0\geq 1$ such that 
$s_{J^{(l_0+p)}}^*\Omega=s_{J^{(l_0)}}^*\Omega$.
Let $L:=J^{(l_0)}$ and 
$\Omega':=s_{L}^*\Omega$.
Then 
$s_{J'}^*\Omega'=\Omega'$ for some $J'$.
This implies that  $\omega_L$ is a sub-Cuntz state by $z=e_{J'}$.
If $\omega$ is pure,
then the statement holds by the construction of $\omega_L$.
\qedh

%\noindent
%We call such a state $\omega$ in Fact \ref{fact:vector} {\it permutative}.
%Fact \ref{fact:vector} shows that 
%any permutative finitely correlated state is quasi-equivalent to 
%a sub-Cuntz state. 

%%%%%%%%%%%%%%%%%%%%%%%%%%%%%%%%%%%%
%
% subsection 4.3
%
\ssft{Infinitely correlated states as non-sub-Cuntz states}
\label{subsection:fourththree}
From Lemma \ref{lem:omegasub}(i),
any infinitely correlated state is not a sub-Cuntz state on $\con$
when $n<\infty$.
In this subsection, 
we show examples of infinitely correlated states.
%
% Example 4.8
%
\begin{ex}
\label{ex:permurep}
(Infinitely correlated state associated with a permutative representation)
{\rm
Let ${\Bbb N}:=\{1,2,\ldots\}$
and let $\{e_{k,m}:(k,m)\in {\Bbb N}\times {\Bbb Z}\}$
denote the standard basis of $\ell^2({\Bbb N}\times {\Bbb Z})$.
For $2\leq n<\infty$,
define a representation $\pi$ of $\con$ 
on $\ell^2({\Bbb N}\times {\Bbb Z})$ by
%
% Equation 4.11
%
\begin{equation}
\label{eqn:pisie}
\pi(s_i)e_{k,m}:=e_{n(k-1)+i,m+1}\quad((k,m)\in {\Bbb N}\times {\Bbb Z},\,
i=1,\ldots,n).
\end{equation}
By definition,
$(\ell^2({\Bbb N}\times {\Bbb Z}),\pi)$ is permutative, and 
$\pi(s_1^m)^*e_{1,0}=e_{1,-m}$ for any $m\geq 1$.
Hence $\dim{\rm Lin}\langle\{\pi(s_J)^*e_{1,0}:J\}\rangle=\infty$.
Therefore 
the state $\omega:=\langle e_{1,0}|\pi(\cdot)e_{1,0}\rangle$
is infinitely correlated.
}
\end{ex}
%
% Example 4.9
%
\begin{ex}
\label{ex:quasifree}
{\rm
(Quasi-free states)
We show that any quasi-free state on $\con$ is infinitely correlated.
Let $\Lambda_n:
=\{(a_{1},\ldots,a_n)\in{\Bbb R}^n: a_i>0 \mbox{ for all }i,\, 
a_1+\cdots+a_n=1\}$.
For $a\in\Lambda_n$, define
the state $\rho_{a}$ on $\con$ by
%
% Equation 4.12
%
\begin{equation}
\label{eqn:rhoasj}
\rho_{a}(s_{J}s_K^*):=a_J\delta_{JK}\quad
(J,K\in{\cal I}:=\bigcup_{m\geq 0}\{1,\ldots,n\}^{m})
\end{equation}
where $a_J:=a_{j_1}\cdots a_{j_m}$ for $J=(j_1,\ldots,j_m)$
and $a_{\emptyset}:=1$.
The state $\rho_{a}$ is called the {\it quasi-free state} on $\con$ by $a$
\cite{ACE,Evans}.
It is known that the GNS representation by $\rho_a$
is a type III factor representation;
$\rho_{a} \sim \rho_b$ if and only if $a=b$; 
$\rho_a\ptimes \rho_b=\rho_{a\boxtimes b}$
\cite{Izumi, TS11,TS07}.
%Hence $\Lambda_n$ is the set of all complete invariants
%of quasi-free states on $\con$.

Fix $a\in\Lambda_n$ and 
let $({\cal H},\pi,\Omega)$ denote 
the GNS representation by $\rho_{a}$.
For $J\in {\cal I}$,
define 
$v_J:=a_J^{-1/2}\pi(s_J)^*\Omega\in {\cal H}$.
From (\ref{eqn:rhoasj}),
we see that $\{v_J:J\in {\cal I}\}$ is an orthonormal system in ${\cal H}$.
Therefore
$\dim {\rm Lin}\langle\{\pi(s_J)^*\Omega:J\in {\cal I}\}\rangle=
\dim {\rm Lin}\langle\{v_J:J\in {\cal I}\}\rangle
=\#{\cal I}=\infty$.
Hence $\rho_a$ is infinitely correlated.
}
\end{ex}

%\ww {\bf Acknowledgment:}
%The author would like to express his sincere thanks to ***
%for his interest in this topic and raising the above question.

%%%%%%%%%%%%%%%%%%%%%%%%%%%%%%%
\appendix \section*{Appendix}
%%%%%%%%%%%%%%%%%%%%
%
% Appendix A
%
\sftt{%Applied 
Combinatorics on words
in tensor algebra
%Periodicity in tensor semigroup
}
\label{section:appone}
In this section,
we prove the freeness of some semigroup
with uncountable rank
associated with a tensor algebra.
By using this fact and known results about free semigroups,
we derive crucial lemmas for main theorems.
%%%%%%%%%%%%%%%%%%%%%%%%%%
%
% subsection A.1
%
\ssft{Cancellation law and equidivisibility of tensor product}
\label{subsection:apponeone}
We show the cancellation law (\cite{Howie}, 2.6.1) and 
equidivisibility (\cite{Howie}, $\S$ 7.1) of tensor product.
Let $U(1):=\{c\in {\Bbb C}:|c|=1\}$.
The following is very elementary, but
mistakes are often found in literature.
%
% Lemma A.1
%
\begin{lem}
\label{lem:trivial}
If $z$ and $y$ are unit vectors in a Hilbert space, then
$\langle z|y\rangle =1$ if and only if $z=y$.
\end{lem}
%
% Proof
% 
\pr
If $\langle z|y\rangle =1$,
then $1=|\langle z|y\rangle|\leq \|z\|\cdot \|y\|=1$.
Since $|\langle z|y\rangle|= \|z\|\cdot \|y\|$,
$z=cy$ for some $c\in U(1)$.
From this,
$1=\langle z|y\rangle =\langle cy|y\rangle =\overline{c}$.
Hence $z=y$.
The inverse direction is trivial.
\qedh
%
%
% Lemma A.2
%
\begin{lem}
\label{lem:cancel}
Let $({\cal V}_{n,m})_1$ be as in $\S$\ref{subsection:firsttwo}.
\begin{enumerate}
%(i)
\item
(Cancellation law)
Assume that $x,x'\in ({\cal V}_{n,m})_1$
and $y,z\in \bigcup_{a\geq 1}({\cal V}_{n,a})_1$ satisfy
$x\otimes y=x'\otimes z$ or $y\otimes x=z\otimes x'$.
Then $x=cx'$ and
$y=\overline{c}z$ for some $c\in U(1)$.
In addition,
if $x=x'$, then $y=z$.
%(ii)
\item
(Equidivisibility)
Assume 
$x\in ({\cal V}_{n,m})_1$,
$y\in ({\cal V}_{n,l})_1$
and $m>l$.
\begin{enumerate}
%(a)
\item
If $x\otimes w=y\otimes z$,
then there exists $x'\in ({\cal V}_{n,m-l})_1$
 such that
$x=y\otimes x'$.
%(b)
\item
If $w\otimes x=z\otimes y$,
then there exists $x''\in ({\cal V}_{n,m-l})_1$
 such that
$x=x''\otimes y$.
\end{enumerate}
\end{enumerate}
\end{lem}
%
% Proof
%
\pr
(i) By using Lemma \ref{lem:trivial},
the statement can be verified.

\noindent
(ii)
When $x=\sum x_Je_J$ and $y=\sum y_K e_K$,
define
%
% Equation A.1
%
\begin{equation}
\label{eqn:xdashb}
x':=\sum_{|J_1|=l}\,\sum_{|J_2|=m-l}
x_{J_1J_2}\,\overline{y_{J_1}}\,e_{J_2}\in {\cal V}_{n,m-l}.
\end{equation}
Since $x'=\sum_{|K|=m-l}\langle y\otimes e_K|z\rangle e_K$,
%
% Equation A.2
%
\begin{equation}
\label{eqn:btwo}
\|x'\|^2
=\sum_{K}|\langle y\otimes e_K|z\rangle |^2
\leq
\sum_{i}|\langle v_i|z\rangle |^2=\|z\|^2=1
\end{equation}
where  $\{v_i\}$ is an orthonormal basis of ${\cal V}_{n,m}$
such that $\{y\otimes e_K:|K|=m-l\}\subset \{v_i\}$.
Then we can verify 
$\langle x\otimes w|y\otimes z\rangle=
\langle x'\otimes w|z\rangle$.
%From (\ref{eqn:btwo}), $x'\in {\cal V}_{n,m-1}$ is well defined.
From this,
%
% Equation A.3
%
\begin{equation}
\label{eqn:arrayrl}
\begin{array}{rl}
1=&\langle x\otimes w|y\otimes z\rangle\\
= &\langle x'\otimes w|z\rangle\\
= &\langle y\otimes  x'\otimes w|y\otimes z\rangle\\
=& \langle y\otimes  x'\otimes w|x\otimes w\rangle\\
= &\langle y\otimes  x'|x\rangle.\\
\end{array}
\end{equation}
From this and (\ref{eqn:btwo}),
$1=\langle y\otimes  x'|x\rangle \leq 
\|y\otimes  x'\|\cdot \|x\|\leq 1$.
This implies
$y\otimes x'=c x$ for some $c\in {\Bbb C}$.
By substituting this into (\ref{eqn:arrayrl}),
$1=\langle y\otimes  x'|x\rangle=\langle cx|x\rangle=\overline{c}$.
Hence $y\otimes x'=x$ and $\|x'\|=1$.
Hence (a) is proved.
As the same token,
 (b) is verified.
\qedh

%%%%%%%%%%%%%%%%%%%%%%%%%%%
%
% subsection A.2
%
\ssft{Projective homogeneous tensor semigroup is free}
\label{subsection:apponetwo}
Let $V:={\cal V}_{n,1}$ and 
we identify ${\cal V}_{n,m}$ with $V^{\otimes m}$ for $m\geq 1$,
and let $V^{\otimes 0}:={\Bbb C}$.
By forgetting the addition of 
the tensor algebra $T(V):=\bigoplus_{m\geq 0}V^{\otimes m}$
over $V$, 
%\cite{Lang},
% p633
$T(V)$ is regarded as a semigroup with respect to the tensor product $\otimes$.
Furthermore, its projective space 
$PT(V):=(T(V)\setminus \{0\})/{\Bbb C}^{\times}$
is a semigroup
with respect to the product $[x][y]:=[x\otimes y]$
for $x,y\in T(V)\setminus \{0\}$
where $[x]:=\{cx: c\in {\Bbb C}^{\times}\}\in PT(V)$.
For any subsemigroup $S$ of $T(V)$,
$PS:=(S\setminus \{0\})/{\Bbb C}^{\times}$
is a subsemigroup of $PT(V)$.
Especially,
we consider the following subsemigroup $G$ of $(T(V),\otimes)$
and its projective semigroup $PG$:
%
% Equation A.4
%
\begin{equation}
\label{eqn:gpg}
G:=\bigcup_{m\geq 1}(V^{\otimes m})_1,
\end{equation}
that is, $G$ is the subsemigroup of all homogeneous unit vectors in $T(V)$
except vectors in $V^{\otimes 0}$.

A semigroup $S$ is said to be {\it free} if 
there exists a nonempty subset $B$ of $S$
such that $B$ generates $S$, and 
for any semigroup $S'$ and any map $f$ from $B$ to $S'$,
there exists a homomorphism $\hat{f}$ from $S$ to $S'$
such that $\hat{f}|_B=f$
%where $\iota$ denotes the inclusion map of $B$ into $S$
\cite{Howie}.
In this case, $S$ is called the {\it free semigroup over $B$}
and $\#B$ is called the {\it rank} of $S$.
A free semigroup is defined uniquely up to an isomorphism 
by the rank. 
%
% Lemma A.3
%
\begin{lem}
\label{lem:freeequiv}
(\cite{CP02}, %p116,  
Theorem 9.1)
A semigroup $S$ is free if and only if 
there exists a nonempty subset $B$ of $S$
such that  every element of $S$
has a unique expression as a product of elements of $B$.
\end{lem}
%
% Proposition A.4
%
\begin{prop}
\label{prop:free}
For $G$ in (\ref{eqn:gpg}), its projective semigroup $PG$ is free.
\end{prop}
%
% Proof
%
\pr
Let ${\goth I}_n$ be as in Corollary \ref{cor:sc}.
From Lemma \ref{lem:freeequiv},
it is sufficient to show that 
every element of $PG$ can be expressed uniquely as a
product of elements of $B:=P{\goth I}_n$.

Let $X\in PG$.
By definition,  $X\in P({\cal V}_{n,m})_1$ for some $m\geq 1$.
Hence $X=[x]$ for some $x\in ({\cal V}_{n,m})_1$.
If $x\in {\goth I}_n$, then $X\in B$.
Hence $X$ is uniquely written as an element of $B$.
Assume $x\not \in {\goth I}_n$.
By definition,  $x=z_1\otimes z_2$ for some $z_1,z_2\in 
\bigcup_{a\geq 1}({\cal V}_{n,a})_1$.
When  $y\in ({\cal V}_{n,m})_1$, define $|y|:=m$.
Then $1\leq |z_1|,|z_2|<m=|x|$.
By decomposing $x$ repeatedly,
we can obtain a finest decomposition $x=x_1\otimes \cdots\otimes x_l$.
Then $x_i\in {\goth I}_n$ for  $i=1,\ldots,l$ and $l\leq m$.
Hence $X=[x]$ always has a decomposition 
$[x_1]\cdots [x_l]$ for $[x_i]\in B$ for all $i$.
Assume that $X=X_1'\cdots X_k'$
and $X_j'=[x_j']$ for $x_j'\in {\goth I}_n$ for all $j$.
Then $x_1\otimes \cdots \otimes x_l=cx_1'\otimes \cdots \otimes x_k'$ 
for some $c\in U(1)$.
From Lemma \ref{lem:cancel}(ii), $|x_1|=|x_1'|$
because $x_1,x_1'\in {\goth I}_n$.
From this and Lemma \ref{lem:cancel}(i),
$x_{1}=c_1x_1'$ for some $c_1\in U(1)$.
This implies $X_1=X_1'$.
By the mathematical induction, we can verify that $X_i=X_i'$ 
for all $i$ and $l=k$.
Therefore $X$ has a unique expression as a product of elements of $B$.
\qedh

\noindent
Remark that $T(V)$ is the free ${\Bbb C}$-algebra 
over the set $\{1,\ldots,n\}$ when $\dim V=n$ \cite{Lang}.
% p633
On the other hand, $PG$ is the free semigroup
over the uncountable set $P{\goth I}_n$,
that is, the rank of $PG$ is uncountable.
%
% Proposition A.5
%
\begin{prop}
\label{prop:propertiesfree}
Let $B^+$ denote the free semigroup over a set $B$.
\begin{enumerate}
%(i)
\item
(\cite{Howie}, Proposition 7.1.5) Let $\#B\geq 2$, and let $u,v\in B^+$.
Then $uv=vu$ if and only if $u$ and $v$ are powers of the
same element $w\in B^+$.
%(ii)
\item
(\cite{Howie},  Proposition 7.1.6) 
If $u,v\in B^+$ satisfy $u^m=v^n$ for some $m,n\geq 1$,
then $u$ and $v$ are powers of the same element $w\in B^+$.
%$u$ and $v$ are both expressible as powers of some $w\in B^+$.
%(iii)
\item
(\cite{Lothaire}, Proposition 1.3.3)
Recall the definition of conjugacy in Theorem \ref{Thm:mainthree}(ii).
Let $x,y \in B^n:=\{b_1\cdots b_n:b_i\in B,\,i=1,\ldots,n\}$ 
and $s,t$ be nonperiodic such that
$x=s^p$ and $y=t^q$.
Then $x$ and $y$ are conjugate if and only if 
$s$ and $t$ are conjugate. 
\end{enumerate}
\end{prop}

%
% Corollary A.6
%
\begin{cor}
\label{cor:projective}
Let  $x\in ({\cal V}_{n,m})_1$
and $y\in ({\cal V}_{n,l})_1$.
\begin{enumerate}
%(i)
\item
Assume 
%
% Equation A.5
%
\begin{equation}
\label{eqn:commuteb}
y\otimes x=c x\otimes y
\end{equation}
for some $c\in U(1)$.
Then there exists 
$w\in ({\cal V}_{n,a})_1$
such that $x=\gamma_1 w^{\otimes f_1}$ and $y=\gamma_2 w^{\otimes f_2}$
for some
$f_1,f_2\geq 1$ and
$\gamma_1,\gamma_2 \in U(1)$.
Especially, $x\otimes y$ is periodic and $c=1$.
%(ii)
\item
Assume that there exist two integers $\alpha,\beta\geq 1$ such that 
%
% Equation A.6
%
\begin{equation}
\label{eqn:powersb}
x^{\otimes \alpha}=y^{\otimes \beta}.
\end{equation}
Then there exists $w$ such that $x=\gamma_1w^{\otimes k_1}$
and $y=\gamma_2w^{\otimes k_2}$
for some $k_1, k_2\geq 1$ and $\gamma_1,\gamma_2\in U(1)$.
Especially, 
if $m>l$, then $x$ is periodic.
If $x$ is nonperiodic,
then $y=cx^{\otimes d}$
for some $d\geq 1$ and $c\in U(1)$.
%(iii)
\item
Assume that $m>l$ and there exist $z_1,z_2$ such that 
%
% Equation A.7
%
\begin{equation}
\label{eqn:modb}
z^{\otimes \alpha}=z_1\otimes z_2,\quad 
y^{\otimes \beta}=z_2\otimes z_1
\end{equation}
for some $\alpha,\beta$.
Then $z$ is periodic.
\end{enumerate}
\end{cor}
%
% Proof
%
\pr
From Proposition \ref{prop:free} and its proof,
Proposition \ref{prop:propertiesfree}
can be applied to the pair $(B^+,B)=(PG,P{\goth I}_n)$.

\noindent
(i)
From (\ref{eqn:commuteb}),
$[y][x]=[x][y]$ in $PG$.
From Proposition \ref{prop:propertiesfree}(i),
$[x],[y]\in \{W^p:p\geq 1\}$ for some $W\in PG$.
Since $W=[w]$ for some $w\in G$,
we obtain the statement.

\noindent
(ii)
From (\ref{eqn:powersb}),
$[x]^{\alpha}=[y]^{\beta}$ in $PG$.
From Proposition \ref{prop:propertiesfree}(ii),
the statement holds.

\noindent
(iii)
Assume that $z=u^{\otimes p}$
and $y=v^{\otimes q}$
for some nonperiodic elements $u$ and $v$.
From Proposition \ref{prop:propertiesfree}(iii),
$[u]$ and $[v]$ are conjugate in $PG$.
This implies $u,v\in ({\cal V}_{n,k})_1$ for some  $k\geq 1$.
Hence
$z=u^{\otimes p}\in ({\cal V}_{n,kp})_1$
and 
$y=v^{\otimes q}\in ({\cal V}_{n,kq})_1$.
Therefore $m=kp$ and $l=kq$.
Since $m>l$, $p>q\geq 1$.
Therefore $z$ is periodic.
\qedh

%%%%%%%%%%%%%%%%%%%%%%%%%
%
% Appendix B
%
\sftt{Proofs of properties of Cuntz states}
\label{section:apptwo}
In this section, we prove well-known basic properties of Cuntz states
on $\con$ ($2\leq n\leq \infty$) \cite{BJ}. 
Since both Fact \ref{fact:existence} and Theorem \ref{Thm:maintwo}(ii)
are proved by using properties of Cuntz states,
we do not use results of sub-Cuntz states here.
Let ${\Bbb N}:=\{1,2,3,\ldots\}$ and 
let $s_1,s_2,\ldots$ denote the Cuntz generators of $\con$.
Define ${\cal I}:=\bigcup_{a\geq 0}{\cal I}_1^a$ where
${\cal I}_1:=\{1,\ldots,n\}$ when $n<\infty$, and
${\cal I}_1:={\Bbb N}$ when $n=\infty$.
Define $\goh:=\ell^2({\cal I}_1)$
and  $\goh_1:=\{z\in \goh: \|z\|=1\}$.
Here we identify $\goh$ with the set of all complex sequences
$(z_i)_{i\in {\cal I}_1}$ such that $\sum_{i}|z_i|^2<\infty$.
%
% Theorem B.1
%
\begin{Thm}
\label{Thm:cuntzstate}
Fix $2\leq n\leq \infty$.
\begin{enumerate}
%(i)
\item
There exists a unique state $\omega_1$ on $\con$ such that $\omega_1(s_1)=1$.
In this case, $\omega_1$  is pure and $\omega_1(s_i)=0$ when $i\ne 1$.
%(ii)
\item
For any $z\in \goh_1$,
a Cuntz state on $\con$ by $z$
 exists uniquely and is pure.
%(iii)
\item
For $z\in \goh_1$,
let $\omega_z$ denote the Cuntz state by $z$.
Then
$\omega_z\sim \omega_y$ if and only if $z=y$.
%(iv)
\item
For any $J,K$,
$\omega_z(s_Js_K^*)=\overline{z_J}\,z_K$.
\end{enumerate}
\end{Thm}
%
% Proof
%
\pr
(i)
Let $({\cal H},\pi,\Omega)$ denote the GNS representation by $\omega_1$.
Then we see  $\pi(s_i)^*\Omega=\delta_{i 1}\Omega$ for all $i$.
This implies  that 
$\omega_1(s_{J}s_K^*)=1$ when $J,K\in W
:=\{\emptyset, (1),(11),(111),\ldots\}\subset {\cal I}$,
and $\omega_1(s_{J}s_K^*)=0$ otherwise.
Therefore the uniqueness of $\omega_1$ holds.

We prove the existence and purity as follows.
Let $\{e_k:k\in {\Bbb N}\}$ denote the standard basis of $\ltn$.

Assume $n<\infty$.
Define the action of $\con$ on $\ltn$ by
%
% Equation B.1
%
\begin{equation}
\label{eqn:standardc}
s_ie_{k}:=e_{n(k-1)+i}\quad(i=1,\ldots,n,\,k\in {\Bbb N}).
\end{equation}
Since $s_1e_1=e_1$, $\omega_1:=\langle e_1|(\cdot)e_1\rangle$ satisfies
$\omega_1(s_1)=1$.
Therefore the existence is proved.
Next we prove the irreducibility of the action (\ref{eqn:standardc}).
Remark that any $k\in {\Bbb N}$
is uniquely written as $n(k'-1)+i$ for some 
$i=1,\ldots,n$ and $k'\in {\Bbb N}$.
Hence we see $\{e_k:k\in{\Bbb N}\}=\{s_Je_1:J\in {\cal I}\}$.
From this and (\ref{eqn:standardc}),
$e_1$ is a cyclic vector of $\ltn$.
Let $v=\sum_{m\geq 1}v_m e_m\in\ltn$, $v\ne 0$.
Define  $m_0\:={\rm min}\{m\in{\Bbb N}:v_m\ne 0\}$.
Then there exists $J_0\in{\cal I}_1^k$ for some $k\geq 1$
such that $e_{m_0}=s_{J_0}e_1$.
Hence $\langle e_1|s_{J_0}^*v \rangle =v_{m_0}\ne 0$.
Therefore we can construct $v'\in \con v$
such that $v'=e_1+v''$ for some
$v''\in \ltn$, $\langle e_1|v''\rangle=0$.
Then we can verify 
$\|(s_1^{*})^lv'-e_1\|\to 0$ for $l\to \infty$.
Therefore $e_1\in \overline{\con v}$.
This implies that  any non-zero invariant closed subspace of $\ltn$ coincides with
$\ltn$.
Therefore the action in (\ref{eqn:standardc})
is irreducible. Hence $\omega_1$ is pure.

Assume $n=\infty$.
Define the action of $\coni$ on $\ltn$ by
%
% Equation B.2
%
\begin{equation}
\label{eqn:siek}
s_{i}e_k:=e_{2^{i-1}(2k-1)}\quad(i,k\geq 1).
\end{equation}
Then $s_1e_1=e_1$ and $\{s_Je_k:J\}=\{e_m:m\in{\Bbb N}\}$.
Therefore $e_1$ is a cyclic unit vector of $\ltn$ and 
$\omega_1:=\langle e_1|(\cdot)e_1\rangle$
satisfies $\omega_1(s_1)=1$.
As the same token,
we can prove that the action in (\ref{eqn:siek}) is irreducible.
Hence $\omega_1$ is pure.

\noindent
(ii)
Let $U(\goh)$ denote the unitary group on $\goh$.
Let $\{e_i\}$ denote the standard basis of $\goh$.
For $z\in\goh_1$,
let $g=(g_{ij})\in U(\goh)$ such that $ge_1=z$
where $g_{ij}:=\langle e_i|ge_j\rangle$.
Then $g_{j1}=z_j$ for all $j$.
For $\omega_1$ in (i),
define 
%
% Equation B.3
%
\begin{equation}
\label{eqn:omegadash}
\omega':=\omega_1\circ \alpha_{g^*}
\end{equation}
where $\alpha$ is as in (\ref{eqn:alphagsi}),
which can be also well defined when $n=\infty$.
By (\ref{eqn:omegadash}),
$\omega'$ is pure and 
we can verify $\omega'(s_j)=\overline{z_j}$ for all $j$
where we use $\omega_1(s_i)=0$ when $i\ne 1$.
Hence $\omega'$ is a Cuntz state by $z$.
Therefore the existence is proved.

If $\omega''$ is a Cuntz state by $z$,
then we can verify that $(\omega''\circ \alpha_g)(s_1)=1$
for $g$ in (\ref{eqn:omegadash}).
This implies 
that $\omega''\circ \alpha_g=\omega_1$ in (i) and 
$\omega''=\omega_1\circ \alpha_{g^*}=\omega'$.
Hence the uniqueness
of the Cuntz state by $z$ is proved.

\noindent
(iii)
Assume $\omega_z\sim \omega_y$.
Then there exists 
an action of $\con$  on a Hilbert space ${\cal H}$
with two cyclic unit vectors $\Omega_z$ and $\Omega_y$
such that 
$s(z)\Omega_z=\Omega_z$
and 
$s(y)\Omega_y=\Omega_y$
from Lemma \ref{lem:asse} (the case of $n=\infty$ also holds)
where $s(z):=\sum z_js_j\in \con$.
Then
$\langle\Omega_z|s_J\Omega_y\rangle
=\overline{z_J}\langle\Omega_z|\Omega_y\rangle$ for any 
$J\in{\cal I}_1^k$.
Since $\{s_J\Omega_y:J\}$ spans ${\cal H}$,
$\langle\Omega_z|\Omega_y\rangle\ne 0$
because $\Omega_z\ne 0$.
On the other hand,
%
% Equation B.4
%
\begin{equation}
\label{eqn:sysx}
\langle\Omega_z|\Omega_y\rangle
=
\langle s(z)\Omega_z|s(y)\Omega_y\rangle
=
\langle \Omega_z|s(z)^*s(y)\Omega_y\rangle
=
\langle z|y\rangle 
\langle\Omega_z|\Omega_y\rangle.
\end{equation}
Since $\langle\Omega_z|\Omega_y\rangle\ne 0$,
$\langle z|y\rangle$ must be $1$. This implies $z=y$
from Lemma \ref{lem:trivial}.
The inverse direction is trivial.

\noindent
(iv)
By definition, the statement is verified.
\qedh
%%%%%%%%%%%%%%%%%%%%%%%
%
% Appendix C
%
\sftt{Sub-Cuntz states on $\coni$}
\label{section:appthree}
In this section, we generalize sub-Cuntz states on $\con$ ($n<\infty$) to $\coni$.
Except some parts,
main theorems and properties of the state parametrization
hold like the case of $n<\infty$.
Hence we list different points and some remarks for the case of $\coni$.
%%%%%%%%%%%%%%%%%%%%%
%%%%%%%%%%%%%%%%
%
%  subsection C.1
%
\ssft{Definition and parametrization}
\label{subsection:appthreeone}
Let ${\Bbb N}:=\{1,2,\ldots\}$ and 
let $\coni$ denote the {\it Cuntz algebra} \cite{C}, that is, 
a C$^{*}$-algebra which is universally generated 
by $\{s_{i}:i\in {\Bbb N}\}$ satisfying
%
% Equation C.1
%
\begin{equation}
\label{eqn:coni}
s_{i}^{*}s_{j}=\delta_{ij}I\quad(i,j\in {\Bbb N}),\quad 
\sum_{i=1}^{k}s_{i}s_{i}^{*}\leq I\quad \mbox{for any }k\in {\Bbb N}.
\end{equation}
For a unit vector $z\in \ltn$,
$\omega_z$ is a {\it Cuntz state} on $\coni$ by $z$ if 
$\omega_z$ is a state on $\coni$ which 
satisfies $\omega_z(s_{i})=\overline{z}_i$ for all $i$.
Then $\omega_z$ exists uniquely 
and is pure for any $z$;
$\omega_z\sim \omega_y$
if and only if $z=y$ (see Appendix \ref{section:apptwo}).
%
% Theorem C.1
%
\begin{Thm}
\label{Thm:conimain}
For $m\geq 1$ and a unit vector $z=\sum z_Je_J\in \ell^2({\Bbb N}^m)$,
there exists a state $\omega$ on $\coni$ which satisfies
$\omega(s_{J})=\overline{z_J}$ for all $J\in{\Bbb N}^{m}$.
Such $\omega$ is called a sub-Cuntz state by $z$ of order $m$.
\end{Thm}
%
% Proof
%
\pr
Fix a bijection $f:{\Bbb N}\cong {\Bbb N}^m$ and 
define the endomorphism $\hat{f}$ of $\coni$ by
$s_{i}\mapsto \hat{f}(s_{i}):=s_{f(i)}$ for each $i\in {\Bbb N}$
where 
$s_{f(i)}:=s_{j_1}\cdots s_{j_m}$ when $f(i)=(j_1,\ldots,j_m)$.
Then $\hat{f}(\coni)\cong \coni$
and $\hat{f}(s_{i})$'s are Cuntz generators of  $\hat{f}(\coni)$.
Then,
for a unit vector $z\in \ell^2({\Bbb N}^m)$,
$\omega$ is a sub-Cuntz state on $\coni$ by $z$ 
if and only if 
$\omega$ is an extension of the Cuntz state $\omega_{\hat{z}}$
on $\hat{f}(\coni)$ by $\hat{z}$
to $\coni$ where 
$\hat{z}:=(z_{f(i)})_{i\in{\Bbb N}}\in \ltn$.
Since an extension of $\omega_{\hat{z}}$ to $\coni$ always exists
from Proposition \ref{prop:stateext}(i),
the statement holds.
\qedh

Let ${\cal V}_{\infty,m}:=\ell^2({\Bbb N}^m)$ with 
the standard basis $\{e_J:J\in {\Bbb N}^m\}$.
We identify ${\cal V}_{\infty,m}$ with $\ell^2({\Bbb N})^{\otimes m}$.
Then the periodicity, decomposability and equivalence 
of parameters are defined as same as  the case of $n<\infty$.

Let 
$\omega$ be  a state on $\coni$ with the GNS representation
$({\cal H},\pi,\Omega)$ and $z=\sum z_Je_{J}\in 
(\ell^2({\Bbb N}^m))_1$.
In a similar fashion  with Theorem \ref{Thm:bj},
we see that
the following are equivalent:
(i) $\omega$ is a sub-Cuntz state by $z$.
(ii)
$\sum z_J\pi(s_J)\Omega=\Omega$.
(iii) $\pi(s_J)^*\Omega=z_J\Omega$ for all $J$.
Remark that $s(z)=\sum z_Js_J$ is well defined in $\coni$
for any $z$.
Hence the l.h.s in (ii) is well defined.
%%%%%%%%%%%%%%%%%%%%%%%%
%
% subsection C.2
% 
\ssft{Main theorems and naturalities of parametrization}
\label{subsection:appthreetwo}
Statements in main theorems are almost same with the case of $n<\infty$.
Let $U(\infty)$ denote the group
of all unitaries in ${\cal B}(\ltn)$.
Then the state parametrization $z\mapsto 
\tilde{\omega}_{z}$
satisfies the $U(\infty)$-covariance.

For the $\varphi$-tensor product of states on $\coni$,
the following new definitions are necessary.
For $2\leq n<\infty$,
let $\{s^{(\infty)}_k\}$ 
and $\{s_i^{(n)}\}$ denote
Cuntz generators of $\coni$ and $\con$, respectively.
Define the embedding $\varphi_{\infty,n}$ of $\coni$ into $\coni\otimes \con$ by
$\varphi_{\infty,n}(s^{(\infty)}_{n(k-1)+i}):=s^{(\infty)}_k\otimes s_i^{(n)}$ 
for $k\in{\Bbb N},\, i=1,\ldots,n$.
For $\omega_1\in{\cal S}_{\infty}$
and $\omega_2\in{\cal S}_n$,
define
$\omega_{1}\ptimes \omega_2\in {\cal S}_{\infty}$ by
$\omega_{1}\ptimes \omega_2
:=(\omega_1\otimes \omega_2)\circ \varphi_{\infty,n}$.
For $\omega_1\in{\cal S}_{\infty}$,
$\omega_2\in{\cal S}_n$
and $\omega_3\in{\cal S}_{n'}$,
$(\omega_1\ptimes \omega_2)\ptimes \omega_3
=\omega_1\ptimes(\omega_2\ptimes \omega_3)$ 
(\cite{TS02}, the proof of Theorem 1.2(iv)).
If $\omega_1$ and $\omega_2$ are sub-Cuntz states
on $\coni$ and $\con$,  respectively,
then
we see that 
$\omega_1\ptimes \omega_2$ is a sub-Cuntz state on $\coni$.

For $J=(j_1,\ldots,j_m)\in {\Bbb N}^{m}$
and 
$K=(k_1,\ldots,k_m)\in \{1,\ldots,n\}^{m}$,
define
$J\boxtimes K=(l_1,\ldots,l_m)\in{\Bbb N}^{m}$ by
$l_t:=n(j_t-1)+k_t$ for $t=1,\ldots,m$.
For $z\in {\cal V}_{\infty,m}$
and  $y\in {\cal V}_{n,m}$,
define $z\boxtimes y\in {\cal V}_{\infty,m}$ by
%
% Equation C.2
%
\begin{equation}
\label{eqn:zinf}
z\boxtimes y:=\sum_{J\in{\Bbb N}^{m}}
(z\boxtimes y)_J e_{J},\quad 
(z\boxtimes y)_J:=z_{J'}y_{J''}
\end{equation}
where
$J'\in{\Bbb N}^{m}$
and 
$J''\in\{1,\ldots,n\}^{m}$
are uniquely defined as
$J=J'\boxtimes J''$.
For $z\in {\cal V}_{\infty,m}$
and  $y\in {\cal V}_{n,l}$,
define $z*y\in {\cal V}_{\infty,\alpha m}$ by
$z*y:=z^{\otimes \alpha}\boxtimes y^{\otimes \beta}$
where $\alpha$ and $\beta$ are chosen such that 
$\alpha m=\beta l$ is the least common multiple of $m$ and $l$.
If $z\in {\cal V}_{\infty,m}$
and  $y\in {\cal V}_{n,l}$ are nonperiodic,
then
$z*y$ is also nonperiodic
and 
$\tilde{\omega}_z\ptimes \tilde{\omega}_y
=\tilde{\omega}_{z*y}$.
%%%%%%%%%%%%%%%%%%%%%
%
% subsection C.3
%
\ssft{Infinitely correlated sub-Cuntz states on $\coni$}
\label{subsection:appthreethree}
Lemma \ref{lem:omegasub}(i) does not hold for $\coni$.
We prove that 
a sub-Cuntz state on $\coni$ is not always finitely correlated
by using examples.
%
% Proposition C.2
%
\begin{prop}
\label{prop:infsub}
For a unit vector $x=\sum x_i e_i\in \ltn$,
define
$z:=\sum x_je_{j}\otimes e_j\in 
(\ltn^{\otimes 2})_1$,
$X:=\{i\in{\Bbb N}: x_i\ne 0\}$ and $N:=\#X$.
Assume that $\omega$ is a sub-Cuntz state on $\coni$ by $z$.
Then  $\omega$  is finitely correlated if and only if $N<\infty$.
\end{prop}
%
% Proof
%
\pr
Let $({\cal H},\pi,\Omega)$ denote 
the GNS representation by $\omega$.
By definition,
$\omega (s_is_j)=\delta_{ij}\overline{x_i}$\, for all $i,j$.
This implies
$\sum x_i\pi(s_i^2)\Omega=\Omega$ and 
%
% Equation C.3
%
\begin{equation}
\label{eqn:sumdelta}
\omega(s_is_j^*)
=\sum_{l}
\delta_{il}
\delta_{jl}\,
\overline{x_{i}}\,x_{j}
=\delta_{ij}|x_{i}|^2\quad \mbox{for all }i,j.
\end{equation}
We divide the case of $N\geq 2$ from the case of $N=1$.

Assume $N\geq 2$.
In this case, $z$ is nonperiodic. 
From Theorem \ref{Thm:maintwo}(iii) for $\coni$,
$\langle \Omega|\pi(s_i)^*\Omega\rangle =\omega(s_i^*)=0$ for all $i$.
From this and (\ref{eqn:sumdelta}),
$\{\Omega, |x_{i}|^{-1}\pi(s_i^*)\Omega:i\in X\}_{i\geq 1}$ 
is an orthonormal family in ${\cal H}$.
Define ${\cal I}:=\bigcup_{a\geq 0}{\Bbb N}^{a}$.
Then ${\cal I}$ is a free monoid with respect to the concatenation
\cite{Lothaire}.
Define the subsemigroup $W$ of ${\cal I}$ generated by $\{(ii): i\in{\Bbb N}\}$:
%
% Equation C.4
%
\begin{equation}
\label{eqn:wlang}
W:=\langle \{(ii):i\in{\Bbb N}\}\rangle \subset {\cal I}.
\end{equation}
Then 
%
% Equation C.5
%
\begin{equation}
\label{eqn:pisjinf}
\pi(s_J^*)\Omega= 
\left\{
\begin{array}{ll}
x[J_1]\pi(s_i)^*\Omega \quad  &\mbox{when }J=J_1i,\, J_1\in W,\\
\\
x[J]\Omega \quad  &\mbox{when }J\in W,\\
\\
0 &\mbox{otherwise}
\end{array}
\right.
\end{equation}
where
$x[J]:=x_{i_1}x_{i_2}\cdots x_{i_l}\in{\Bbb C}$
when $J=(i_1i_1\, i_2i_2\cdots\, i_li_l)\in W$. 
From this,
${\rm Lin}\langle\{\pi(s_J)^*\Omega:J\}\rangle =
{\rm Lin}\langle\{\Omega,\pi(s_i)^*\Omega:i\in X\}\rangle$.
Therefore  
%
% Equation C.6
%
\begin{equation}
\label{eqn:kdim}
\dim {\cal K}=
\dim {\rm Lin}\langle\{\Omega,\pi(s_i)^*\Omega:i\in X\}\rangle
=1+\#X=1+N
\end{equation}
where ${\cal K}:=\overline{{\rm Lin}\langle\{\pi(s_J)^*\Omega:J\}\rangle}\subset
{\cal H}$.
From (\ref{eqn:kdim}),
the statement holds except the case of $N=1$.

Assume $N=1$.
It is sufficient to show that $\omega$ is finitely correlated.
By assumption,
there exists $j$ such that $z=x_je_j\otimes e_j$ 
and $|x_j|=1$.
In this case,  we obtain 
$x_j\pi(s_j^2)\Omega=\Omega$.
From this,
%
% Equation C.7
%
\begin{equation}
\label{eqn:square}
\pi(s_{i}s_{k})^*\Omega=\delta_{ij}\delta_{jk}x_{j}\Omega
\quad \mbox{for all }i,k.
\end{equation}
Let $q$ be a quadratic root of $x_j$.
Then  there exists $0\leq a\leq 1$ such that
$\omega=a\omega_++(1-a)\omega_-$
where $\omega_{\pm}$ denotes
the Cuntz state by $\pm q e_j$.
In the proof of  Theorem \ref{Thm:periodic},
there exists a unit vector $(\alpha,\beta)\in {\Bbb R}^2$
such that 
$\Omega=\alpha \Omega_+ +\beta \Omega_-$ and 
$\alpha^2=a$ and $\beta^2=1-a$
where $\Omega_{\pm}$ denotes the 
GNS cyclic vector by $\omega_{\pm}$.
Then we see that 
$\pi(s_i)^*\Omega=q\delta_{ij}( \alpha \Omega_+ - \beta \Omega_-)$
for all $i$.
From this and  (\ref{eqn:square}),
$\pi(s_J)^*\Omega\in{\rm Lin }\langle\{\Omega,
\alpha \Omega_+ - \beta \Omega_-\}\rangle$ for any $J$.
Therefore
$\dim {\cal K}\leq \dim {\rm Lin }\langle\{\Omega,
\alpha \Omega_+ - \beta \Omega_-\}\rangle\leq 2$.
Hence $\omega$ is finitely correlated.
\qedh

\noindent
For example, 
let $x:=\sum 2^{-i/2}e_i\in (\ltn)_1$.
%$\|x\|=\sum_{j=1}^{\infty}2^{-j}=
%\frac{1}{2}\frac{1}{1-1/2}=1$ and 
Then $\omega$ 
associated with $x$ satisfies $N=\infty$ in Proposition \ref{prop:infsub}.

%%%%%%%%%%%%%%%%%%%%%%%%%%%%%%%%%%%%%%
%
% Reference 
%

%
\label{Lastpage}


\begin{thebibliography}{99}
%
% AAAAAA
\bibitem{AK03}M.\ Abe, K.\ Kawamura,
Pseudo-Cuntz algebra and recursive FP ghost system in string theory,
Int.\ J.\ Mod.\ Phys.\ A 18(4)  (2003), 607--625.
%
\bibitem{AK02RR}M.\ Abe, K.\ Kawamura,
Branching laws for endomorphisms of fermions 
and the Cuntz algebra ${\cal O}_{2}$,
J.\ Math.\ Phys.\ 49 (2008), 043501-01--043501-10. 
%
\bibitem{ACE}H.\ Araki, A.L.\ Carey, D.E.\ Evans, 
On $O_{n+1}$, 
J.\ Operator Theory 12 (1984),  247--264.
%
\bibitem{BC}W.R.\ Bergmann, R.\ Conti,  
Induced product representation of extended Cuntz algebras,
Annali di Mathematica 182  (2003), 271--286. 
%
\bibitem{Bhatia}R.\ Bhatia,
Matrix analysis, 
Springer, 1997.
%
\bibitem{BJ1997}O.\ Bratteli, P.E.T.\ Jorgensen,
Endomorphisms of ${\cal B}({\cal H})$ II. 
Finitely correlated states on $\con$,
J.\ Funct.\ Anal.\ 145 (1997), 323--373.
%
\bibitem{BJ}O.\ Bratteli,  P.E.T.\ Jorgensen, 
Iterated function systems and permutation representations 
of the Cuntz algebra,
Mem.\ Amer.\ Math.\ Soc.\ {\rm 139} (1999), 1--89.
%
\bibitem{BJKW}O.\ Bratteli,  P.E.T.\ Jorgensen, 
A.\ Kishimoto,   R.F.\ Werner, Pure states on $\co{d}$, 
J.\ Operator Theory 43(1) (2000), 97--143. 
%
\bibitem{BJO}O.\ Bratteli, P.E.T.\ Jorgensen, V.\ Ostrovsky\u{\i},
Representation theory and numerical AF-invariants.
The representations and centralizers of certain states on ${\cal O}_d$, 
Mem.\ Amer.\ Math.\ Soc.\ 168(797) (2004), 1--178.
%
\bibitem{BJP}O.\ Bratteli, P.E.T.\ Jorgensen, G.L.\ Price, 
Endomorphisms of ${\cal B}({\cal H})$,
in
Quantization, nonlinear partial differential equations, and operator algebra 
(W.\ Arveson, T.\ Branson, and I.\ Segal, eds.), 
Proc.\ Sympos.\ Pure Math.,  59, Amer.\ Math.\ Soc., 1996, pp. 93--138.
%
\bibitem{CP02}A.H.\ Clifford, G.B.\ Preston,
The algebraic theory of semigroup  vol.\ II,
American Mathematical Society, 1967.
%
\bibitem{C}J.\ Cuntz,
Simple $C^*$-algebras generated by isometries,
Commun.\ Math.\ Phys.\ 57  (1977), 173--185.
%
\bibitem{DaPi2}K.R.\ Davidson, D.R.\ Pitts, 
The algebraic structure of non-commutative analytic Toeplitz algebras,
Math.\ Ann.\ {\rm 311}  (1998), 275--303.
%
\bibitem{DaPi3}K.R.\ Davidson, D.R.\ Pitts, 
Invariant subspaces and hyper-reflexivity for free semigroup algebras, 
Proc.\ London Math.\ Soc.\ {\rm 78}  (1999), 401--430.	
%
\bibitem{Dixmier}J.\ Dixmier,
$C^*$-algebras,
North-Holland Publishing Company, 1977.
%
\bibitem{DS2}N.\ Dunford,  J.T.\ Schwartz, 
Linear operators.\ II,
Interscience, New York, 1963.
%
\bibitem{DHJ}D.E.\ Dutkay,  J.\ Haussermann,  P.E.T.\ Jorgensen,
Atomic representations of Cuntz algebras,
arXiv:1311.5265v1.
%
\bibitem{ES}L.\ Eld\'{e}n, B.\ Savas,
A Newton-Grassmann method for computing the best multilinear 
rank-$(r_1,$ $r_2,$ $r_3)$ approximation of a tensor,
SIAM.\ J.\ Matrix Anal. \& Appl., 31(2) (2009), 248--271.
%
\bibitem{Evans}D.E.\ Evans,
On $\con$,
Publ.\ RIMS, Kyoto Univ.\ 16 (1980), 915--927.
%
\bibitem{FL}N.J.\ Fowler, M.\ Laca,
Endomorphisms of ${\cal B}({\cal H})$,
extensions of pure states, and a class of representations of $\con$,
J.\ Operator Theory 40(1) (2000), 113--138.
%
\bibitem{Gabriel}M.J.\ Gabriel,
Cuntz algebra states defined by implementers of endomorphisms of the
CAR algebra,
Canad.\ J.\ Math. 54 (2002), 694--708.
%
\bibitem{Glimm}J.\ Glimm,
Type I C$^{*}$-algebras,
Ann.\ Math.\ 73(3)  (1961), 572--612.
%
\bibitem{GV}G.H.\ Golub, C.F.\ Van,
Matrix computations 3rd ed.,
The Johns Hopkins University Press, 1996.
%
\bibitem{Howie}J.M.\ Howie, 
Fundamentals of semigroup theory,
Oxford Science Publications, 1995.
%
\bibitem{Izumi}M.\ Izumi,   
Subalgebras of infinite C$^{*}$-algebras with 
finite Watatani indices.\ I. Cuntz algebras,
Commun.\ Math.\ Phys.\ 155(1) (1993), 157--182.
%
\bibitem{Jeong1999}E.-C.\ Jeong,
Irreducible representations of the Cuntz algebra $\con$, 
Proc.\ Amer.\ Math.\ Soci.\ 127(12) (1999), 3583--3590.
%
\bibitem{Jeong2005}E.-C.\ Jeong,
Linear functionals on the Cuntz algebra,
Proc.\ Amer.\ Math.\ Soc.\ 134(1) (2005), 99--104.
%
%%%%%%%%%%%%%%%%%%%%%%%%%%%%%%%%%%%
%
\bibitem{GP0123}K.\ Kawamura,
Generalized permutative representations of the Cuntz algebras,
arXiv:math/0505101.
%math.OA/0505101.
%
\bibitem{IWF}K.\ Kawamura,
Extensions of representations of the CAR algebra to
the Cuntz algebra ${\cal O}_2$ ---the Fock and the infinite wedge---,
J.\ Math.\ Phys.\ 46(7)   (2005),  073509-1--073509-12.
%
\bibitem{PFO01}K.\ Kawamura,
The Perron-Frobenius operators, invariant measures 
and representations of the Cuntz-Krieger algebras,
J.\ Math.\ Phys.\ 46(8) (2005), 083514-1--083514-6. 
%
\bibitem{PE01}K.\ Kawamura, 
Polynomial endomorphisms of  the Cuntz algebras arising from
permutations. I ---General theory---,
Lett.\ Math.\ Phys.\ 71 (2005), 149--158. 
%
\bibitem{PE02}K.\ Kawamura, 
Branching laws for polynomial endomorphisms of Cuntz algebras
arising from permutations,
Lett.\ Math.\ Phys.\ 77 (2006), 111--126.
%
\bibitem{TS01}K.\ Kawamura,
A tensor product of representations of Cuntz algebras,
Lett.\ Math.\ Phys.\ 82(1) (2007), 91--104.
%
\bibitem{PE03}K.\ Kawamura, 
Automata computation of branching laws for endomorphisms of 
Cuntz algebras,
Int.\ J.\ Alg.\ Comput.\ 17(7) (2007), 1389--1409.
%
\bibitem{TS02}K.\ Kawamura, 
C$^{*}$-bialgebra defined by the direct sum of Cuntz algebras,
J.\ Algebra 319 (2008), 3935--3959.
%
\bibitem{TS11}K.\ Kawamura,
Classification and realizations of type III factor representations
of Cuntz-Krieger algebras associated with quasi-free states,
Lett.\ Math.\ Phys.\ 87 (2009), 199--207.
%
\bibitem{TS08}K.\ Kawamura, 
Pentagon equation arising from state equations of a C$^{*}$-bialgebra,
Lett.\ Math.\ Phys.\ 93 (2010), 229--241.
%math.OA/0906.2507v1.
%
\bibitem{TS15}K.\ Kawamura,
$R$-matrices and the Yang-Baxter equation 
on GNS representations of C$^{*}$-bialgebras,
Linear Alg.\ Appli.\ 438 (2013), 573--583.
%
\bibitem{TS07}K.\ Kawamura,
Tensor products of type III factor 
representations of Cuntz-Krieger algebras,
Algebr.\ Represent.\ Theor.\ 16(5)  (2013), 1397--1407. 
%
%\bibitem{TS19}K.\ Kawamura, 
%C$^{*}$-bialgebra defined as the direct sum of 
%UHF algebras, Comm.\ Algebra, to appear.
%
\bibitem{CFR02}K.\ Kawamura, Y.\ Hayashi, D.\ Lascu,
Continued fraction expansions and permutative representations 
of the Cuntz algebra ${\cal O}_{\infty}$,
J.\ Number Theory 129 (2009), 3069--3080.
%%%%%%%%%%%%%%%%%%%%%%%%%%%%%%%%%%%%%%%
%
\bibitem{Kobayashi}T.\ Kobayashi,
Theory of discretely decomposable restrictions of unitary representations of semisimple 
Lie groups and some applications, 
(translation of S\={u}gaku 51(4) (1999), 337--356),
Sugaku Expositions 18(1) (2005), 1--37.
%
\bibitem{Laca1993}M.\ Laca, 
Gauge invariant states of ${\cal O}_{\infty}$,
J.\ Operator Theory 30(2) (1993), 381--396.
%
\bibitem{Lang}S.\ Lang, 
Algebra Revised 3rd ed., 
Springer, 2002.
%
\bibitem{LMV}L.D.\ Lathauwer, B.D.\ Moor, J.\ Vandewalle,
A multilinear singular value decomposition,
SIAM\ J.\ Matrix\ Anal.\ Appl., 21 (2000), 1253--1278.
%
\bibitem{Lawson2009}M.V.\ Lawson,
Primitive partial permutation 
representations of the polycyclic monoids
and branching function systems,
Periodica Mathematica Hungarica 58 (2) (2009), 189--207.
%
\bibitem{LS}J.-R.\ Lee, D.-Y.\ Shin,
The positivity of linear functionals on Cuntz algebras
associated to unit vectors,
Proc. Amer. Math. Soc. 132(7) (2004), 2115--2119.
%
\bibitem{Lothaire}M.\ Lothaire, 
Combinatorics on words, 
Cambridge University Press, 1983.
%
\bibitem{Shin}D.-Y.\ Shin,
State extensions of states on 
UHF$_n$ algebra to Cuntz algebra,
Bull.\ Korean Math.\ Soc.\ 39(3) (2002),  471--478.
%
\bibitem{VT}M.A.O.\ Vasilescu, D.\ Terzopoulos,
Multilinear analysis of image ensembles: Tensorfaces. 
In Proc. 7th European Conference on Computer Vision (ECCV'02),
Lecture Notes in Computer Science, 2350 (2002),  447--460.
%
\end{thebibliography}
\end{document}